\documentclass{LMCS}
\usepackage{etex,subfig,comment,xspace}
\usepackage{mathrsfs}
\usepackage{amsmath}
\usepackage{url}

\usepackage{%
  amsfonts%
}

\DeclareMathOperator{\bicolimString}{bicol}
\newcommand{\bicolim}[1][]{\ensuremath{\bicolimString_{#1}}}

\DeclareMathOperator{\colim}{col}

\DeclareMathOperator{\Hom}{Hom}

\DeclareMathOperator{\Ob}{ob}

\newcommand{\spa}[2]{%
  \ensuremath{(#1,#2)}%
}

\renewcommand{\propto}[3][]{%
  \ifthenelse{\equal{#1}{}}%
  {\ensuremath{(%
      #2,\,#3
      )}}%
  {\ensuremath{(%
      #2,\,#1,\,#3
      )}}%
}%

\newcommand{\Defeq}
 {\stackrel{\mathrm{def}}{=}}

\newcommand{\Cat}{\mathbf{Cat}}

\newcommand{\tuple}[1]{\langle {#1} \rangle}

\newcommand{\vco}{%
  \ensuremath{\mathbin{%
      \settowidth{\dimen7}{\mbox{$\circ$}}%
      \makebox[0pt][l]{$\circ$}%
      \makebox[\dimen7]{$\cdot$}%
    }}%
}

\newcommand{\natra}[1]{%
  \ensuremath{\mathfrak{#1}}%
}

\DeclareRobustCommand{\cocone}[1]{%
  \natra{#1}%
}

\newcommand{\diag}[2][]{%
  \ensuremath{\mathcal{#2}}
}

\newcommand{\Gr}{%
  \Gamma%
}

\newcommand{\Span}[1][]{%
  \ensuremath{\mathrm{Span}{(\CAT{#1})}}%
}

\newcommand{\Spc}[2][]{%
  {\ensuremath{\mathrm{Sp}(\CAT{#2})}}%
}

\newcommand{\id}{%
  \ensuremath{\mathrm{i\kern-.02em d}}%
}
\newcommand{\CAT}[1]{%
  \ensuremath{\mathbf{#1}}%
}
\newcommand{\cat}{%
  \CAT%
}

\newcommand{\longCAT}[2][]{%
  \ensuremath{\mathbf{#2}_{#1}}%
}

\usepackage{xcolor,graphicx}
\definecolor{%
  catcolor}{gray}{.45}
\newcommand{\Set}{%
  \longCAT{Set}%
}

\DeclareRobustCommand{\SPAN}[1]{%
  \ensuremath{\Span[{\CAT{#1}}]}%
}

\newcommand{\pbDia}{%
  \raisebox{.4ex}{\rotatebox[origin=c]{180}{\scalebox{1.1}{\ensuremath{%
          {\smaller{{\mathrm{\kern-.1 em}\uparrow\mathrm{\kern-.1 em}}}}{%
            \makebox[0pt][l]{\raisebox{-.6ex}{\scalebox{.7}{\ensuremath{\urcorner}}}}%
          }{\vphantom{\raisebox{-.2ex}{g}}}^{\to}_{\to}{\smaller{{\mathrm{\kern-.1 em}\uparrow\mathrm{\kern-.1 em}}}}%
        }}}}%
}%

\newcommand{\pbDiaR}{\reflectbox{\pbDia}}

\newcommand{\pbSqrR}[4]{%
  \ensuremath{{}_{\smaller{#1}}^{\smaller{#4}}\mkern2mu\pbDiaR\mkern3mu{}^{\smaller{#3}}_{\smaller{#2}}}%
}

\let\oldvisiblespace\textvisiblespace
\renewcommand{\textvisiblespace}{\ensuremath{\text{\oldvisiblespace}}}
\newcommand{\strt}[1][x]{%
  \ifthenelse{\equal{#1}{x}}%
  {\ensuremath{\vphantom{|^-_-}}}%
  {\ifthenelse{\equal{#1}{b}}%
    {\ensuremath{\vphantom{|_-}}}%
    {\ifthenelse{\equal{#1}{t}}%
      {\ensuremath{\vphantom{|^-}}}%
      {\ensuremath{\vphantom{|}}}%
    }%
  }%
}

\let\myoldDelta\Delta
 \renewcommand{\Delta}[1][]{%
   \ifthenelse{\equal{#1}{}}%
   {\myoldDelta}%
   {\ensuremath{\myoldDelta}}
}

\usepackage{ifthen%
}
\usepackage{stmaryrd}
\usepackage{amssymb}

\let\oldrightarrow\to
\renewcommand{\to}[1][]{%
  \ifthenelse{\equal{#1}{}}%
  {\oldrightarrow}%
  {\mathbin{%
      \mathchoice%
      {\scalebox{.8}[1]{$\displaystyle\relbar$}%
        {\raisebox{.23ex}{$\textstyle #1$}}%
        {\displaystyle\shortrightarrow}}%
      {\scalebox{.8}[1]{$\textstyle\relbar$}%
        {\raisebox{.23ex}{$\scriptstyle #1$}}%
        {\textstyle\shortrightarrow}}%
      {\scalebox{.8}[1]{$\scriptstyle\relbar$}%
        {\raisebox{.23ex}{$\scriptscriptstyle #1$}}%
        {\scriptstyle\shortrightarrow}}%
      {\scalebox{.8}[1]{$\scriptscriptstyle\relbar$}%
        {\raisebox{.23ex}{$\scriptscriptstyle #1$}}%
        {\scriptscriptstyle\shortrightarrow}}%
  }}%
}

\let\oldleftarrow\gets
\renewcommand{\gets}[1][]{%
  \ifthenelse{\equal{#1}{}}%
  {\oldleftarrow}%
  {\mathbin{%
      \mathchoice%
      {{\displaystyle\shortleftarrow}%
        {\raisebox{.23ex}{$\textstyle #1$}}%
        \scalebox{.8}[1]{$\displaystyle\relbar$}}%
      {{\textstyle\shortleftarrow}%
        {\raisebox{.23ex}{$\scriptstyle #1$}}%
        \scalebox{.8}[1]{$\textstyle\relbar$}}%
      {{\scriptstyle\shortleftarrow}%
        {\raisebox{.23ex}{$\scriptscriptstyle #1$}}%
        \scalebox{.8}[1]{$\scriptstyle\relbar$}}%
      {{\scriptscriptstyle\shortleftarrow}%
        {\raisebox{.23ex}{$\scriptscriptstyle #1$}}%
        \scalebox{.8}[1]{$\scriptscriptstyle\relbar$}}%
  }}%
}

\let\oldrightarrowPart\rightharpoonup
\renewcommand{\rightharpoonup}[1][]{%
  \ifthenelse{\equal{#1}{}}%
  {\oldrightarrowPart}%
  {\mathbin{\scalebox{.8}[1]{$\relbar$}{\raisebox{.23ex}{$\scriptstyle #1$}}{\oldrightarrowPart}}}
}
\let\oldleftarrowPart\leftharpoonup
\renewcommand{\leftharpoonup}[1][]{%
  \ifthenelse{\equal{#1}{}}%
  {\oldleftarrowPart}%
  {\mathbin{{\oldleftarrowPart}{\raisebox{.23ex}{$\scriptstyle #1$}}\scalebox{.8}[1]{$\relbar$}}}
}

\renewcommand{\smaller}[1]{%
    \mathchoice%
    {{\textstyle#1}}%
    {{\scriptstyle#1}}%
    {{\scriptscriptstyle#1}}%
    {{\scriptscriptstyle#1}}%
}

\usepackage{microtype}
\usepackage[latin2]{inputenc}
\usepackage[all]{xy}
\SelectTips{cm}{10}
\usepackage{enumerate,hyperref}
\usepackage{tikz}
\newcommand{\mnode}[3][]{%
  \ifthenelse{\equal{#1}{}}%
  {\node[inner sep = 2pt] (#2) at (#3) {${#2}$};}%
  {\node[inner sep = 2pt] (#1) at (#3) {${#2}$};}%
}

\usetikzlibrary{topaths%
}

\newcommand{\ardraw}[3][->]{%
  \ardrawL[#1]{#2}{#3}{ }{ };
}

\newcommand{\ardrawL}[5][]{%
    \ifthenelse{\equal{#4}{}}
    {%
      \ifthenelse{\equal{#1}{}}
      {%
        \begin{scope}[inner sep=1.5pt,fill=white]
          \draw[->] (#2) --node[#4,fill]  {$\scriptstyle#5$} (#3);
        \end{scope}%
      }%
      {%
        \begin{scope}[inner sep=1.5pt,fill=white]
          \draw[->]  (#2) #1   node[#4,fill]  {$\scriptstyle#5$} (#3);
        \end{scope}%
      }%
    }
    {%
      \ifthenelse{\equal{#1}{}}
      {%
        \begin{scope}[inner sep=1.5pt]
          \draw[->] (#2) --node[#4]  {$\scriptstyle#5$} (#3);
        \end{scope}
      }%
      {%
        \begin{scope}[inner sep=1.5pt]
          \draw[->] (#2)  #1  node[#4]  {$\scriptstyle#5$} (#3);
        \end{scope}
      }%
    }
}
\newcommand{\ardrawd}[3][--]{%
  {\tikzstyle{every node}=[inner sep=3pt]
    \draw[-,thick,white,double] (#2) #1  (#3);
    \draw[->] (#2) #1  (#3);}}

\newcommand{\ardrawLd}[5][--]{%
  \ardrawd[#1]{#2}{#3}%
  \ardrawL[#1]{#2}{#3}{#4}{#5}%
}

\usepgflibrary{arrows%
}

\newcommand{\pardrawL}[5][--]{{%
    \ifthenelse{\equal{#4}{}}
    {\begin{scope}[inner sep=3pt,fill=white,fill]
    \draw[>=left to,->] (#2) #1 node[#4,fill=white]  {$\scriptstyle#5$} (#3);
    \end{scope}}
    {\begin{scope}[inner sep=3pt]
    \draw[>=left to,->] (#2) #1 node[#4]  {$\scriptstyle#5$} (#3);
    \end{scope}}
}}

\newcommand{\PardrawL}[5][--]{{%
    \ifthenelse{\equal{#4}{}}
    {\begin{scope}[inner sep=3pt,fill=white,fill]
    \draw[>=right to,->] (#2) #1 node[#4,fill=white]  {$\scriptstyle#5$} (#3);
    \end{scope}}
    {\begin{scope}[inner sep=3pt]
    \draw[>=right to,->] (#2) #1 node[#4]  {$\scriptstyle#5$} (#3);
    \end{scope}}
}}

\newcommand{\pardrawd}[3][--]{%
  {   \tikzstyle{every node}=[inner sep=3pt]
    \draw[-,thick,white,double] (#2) #1  (#3);
    \draw[>=left to,->] (#2) #1  (#3);}}

\newcommand{\pardrawLd}[5][--]{%
  \pardrawd[#1]{#2}{#3}%
  \pardrawL[#1]{#2}{#3}{#4}{#5}%
}

\newcommand{\POCS}[4][.09]{%
  \begin{scope}[scale=#1]
  \draw (#2)++([scale=.7]#4)++([scale=.7]#3) 
  ++(#3) ++(#3) -- ++ (#4) -- ++ (#4) ;
  \draw (#2)++([scale=.7]#3) ++([scale=.7]#4) 
  ++(#4) ++(#4) -- ++ (#3)-- ++ (#3) ;
  \end{scope}
}
\newcommand{\POCSd}[4][.09]{%
  \begin{scope}[line width=1.5pt,white,double]%
    \POCS[#1]{#2}{#3}{#4}%
  \end{scope}%
  \POCS[#1]{#2}{#3}{#4}%
}

\newcommand{\POC}[3][.11]{%
  \POCS[#1]{#2}{#3-45:1}{#3+45:1}
}

\newdimen\myArrowsize
\newdimen\myDoubleDistance
\pgfarrowsdeclare{TwoCellTip}{TwoCellTip}
{
  \myArrowsize=0.2pt
  \advance\myArrowsize by .5\pgflinewidth
  \pgfarrowsleftextend{.5\myDoubleDistance}
  \pgfarrowsrightextend{1.35\myDoubleDistance}
}
{
  \myArrowsize=0.4pt
  \pgfsetlinewidth{0.4pt}
  \pgfsetdash{}{0pt} 
  \pgfsetroundjoin   
  \pgfsetroundcap    
  \pgfpathmoveto{\pgfpoint{.2pt}{.5\myDoubleDistance}}
  \pgfpatharc{234}{252}{13\myArrowsize}
  \pgfusepathqstroke
  \pgfpathmoveto{\pgfpoint{.2pt}{.5\myDoubleDistance}}
  \pgfpatharc{234}{207}{13\myArrowsize}
  \pgfusepathqstroke
  \pgfsetlinewidth{0.25pt}
  \pgfpathmoveto{\pgfpoint{.2pt}{.5\myDoubleDistance}}
  \pgfpatharc{234}{253.5}{13\myArrowsize}
  \pgfusepathqstroke
  \pgfpathmoveto{\pgfpoint{.1pt}{.5\myDoubleDistance}}
  \pgfpatharc{234}{206}{13\myArrowsize}
  \pgfusepathqstroke
  \pgfsetlinewidth{0.15pt}
  \pgfpathmoveto{\pgfpoint{.2pt}{.5\myDoubleDistance}}
  \pgfpatharc{234}{260}{10.5\myArrowsize}
  \pgfusepathqstroke
  \pgfsetlinewidth{0.4pt}
  \pgfpathmoveto{\pgfpoint{.2pt}{-.5\myDoubleDistance}}
  \pgfpatharc{360-234}{360-252}{13\myArrowsize}
  \pgfusepathqstroke
  \pgfpathmoveto{\pgfpoint{.2pt}{-.5\myDoubleDistance}}
  \pgfpatharc{360-234}{360-207}{13\myArrowsize}
  \pgfusepathqstroke
  \pgfsetlinewidth{0.25pt}
  \pgfpathmoveto{\pgfpoint{.2pt}{-.5\myDoubleDistance}}
  \pgfpatharc{360-234}{360-253.5}{13\myArrowsize}
  \pgfusepathqstroke
  \pgfpathmoveto{\pgfpoint{.1pt}{-.5\myDoubleDistance}}
  \pgfpatharc{360-234}{360-206}{13\myArrowsize}
  \pgfusepathqstroke
  \pgfsetlinewidth{0.15pt}
  \pgfpathmoveto{\pgfpoint{.2pt}{-.5\myDoubleDistance}}
  \pgfpatharc{360-234}{360-260}{10.5\myArrowsize}
  \pgfusepathqstroke
}

\newcommand{\artwoL}[5][--]{%
  \draw[thin,double distance=.35ex,-TwoCellTip] (#2) #1 node[#4,inner sep=4pt] {$\scriptstyle #5$} (#3);%
}

\usetikzlibrary{backgrounds,fit}
\usepgflibrary{arrows}

\renewcommand{\Defeq}{:=}
\newcommand{\vk}{\textsc{vk}}
\renewcommand{\circledcirc}{\vco}

\newcommand{\ie}{i.e.\xspace}
\newcommand{\bicat}{\mathscr}
\newcommand{\CSpan}[1]{\ensuremath{\mathrm{Span}_{#1}^{\Leftarrow}}}
\renewcommand{\natra}[1]{%
 #1%
}
\renewcommand{\cocone}[1]{%
 #1%
}

\def\doi{7 (1:14) 2011}
\lmcsheading%
{\doi}
{1--22}
{}
{}
{Oct.~19, 2010}
{Mar.~31, 2011}
{}

\begin{document}

\title[]{Being Van Kampen is a universal property\rsuper*}

\author[T.~Heindel]{Tobias Heindel\rsuper a}	
\address{{\lsuper a}Laboratoire d'Informatique de Paris-Nord, 
Universit{\'e} de Paris, France}	
\email{Tobias.Heindel@lipn.univ-paris13.fr}

\author[P.~Soboci\'{n}ski]{Pawe{\l} Soboci\'{n}ski\rsuper b}	
\address{{\lsuper b}DSSE, Electronics and Computer Science, University of Southampton, United Kingdom}	
\email{ps@ecs.soton.ac.uk}

\titlecomment{{\lsuper*}This paper is an extended version of the CALCO~`09 paper ``Van Kampen colimits
as bicolimits in Span''~\cite{Heindel2009}.}

\begin{abstract}
 Colimits that satisfy 
 the Van Kampen condition have interesting exactness properties.
 We show that the
 elementary presentation of the Van Kampen condition
 is actually a characterisation 
 of a universal property in the associated bicategory of spans.
 %
 The main theorem states that Van Kampen cocones are precisely
 those diagrams in a category
 that induce bicolimit diagrams in its
 associated bicategory of spans, 
 provided that the category has pullbacks and enough colimits.
\end{abstract}

\maketitle\vspace{2mm}


\section*{Introduction}


 \noindent Exactness, or in other words,
 the relationship between limits and colimits in various categories of interest
 is a research topic
 with several applications in theoretical computer science, 
 including the solution of recursive domain equations~\cite{Scott1972},   
 semantics of concurrent programming
 languages~\cite{Winskel1995}
 and the study of formal grammars and transformation systems~\cite{Corradini1997}.
 Researchers have identified several classes of categories
 in which certain limits and colimits relate to each other in useful ways;
 extensive categories~\cite{Carboni1993,Lawvere1991} and adhesive 
 categories~\cite{Lack2005} 
 are two relatively recent examples. Going further back, research on 
 toposes and quasitoposes involved elaborate study of their exactness 
 properties~\cite{Johnstone1977,Wyler1991}.

 Extensive categories~\cite{Carboni1993} have coproducts that are 
 ``well-behaved'' with respect to pullbacks; more concretely, 
 they are disjoint and stable under pullback. 
 Extensivity has been 
 used by mathematicians~\cite{Brown1997}
 and computer scientists~\cite{Milius2003} alike. 
 In the presence of products, 
 extensive categories are distributive~\cite{Carboni1993} 
 and thus can be used, for instance, 
 to model circuits~\cite{Walters1991}
 or to give models of specifications~\cite{Jacobs1995}.
 Sets and topological spaces inhabit extensive categories while
 quasitoposes are not, in general, extensive~\cite{Johnstone2002}.

 Adhesive categories~\cite{Lack2004,Lack2005}
 have pushouts along monos that are similarly
 ``well-behaved'' with respect to pullbacks---they are instances of
 Van Kampen squares.
 Adhesivity has been used as a categorical foundation
 for double-pushout graph transformation~\cite{Lack2004,Ehrig2004a}
 and has found several related applications~\cite{Ehrig2004,Sassone2005a}. 
 Toposes are adhesive~\cite{Lack2006} 
 but quasitoposes, in general, are not~\cite{Johnstone2007}. 

 The elementary characterisations of coproducts in extensive categories 
 and pushouts along monos in adhesive categories can be seen 
 as specific instances of a general condition that can be expressed for any
 colimit. Cockett and Guo dubbed the colimits that satisfy this condition  
 \emph{Van Kampen \textsc{(vk)} colimits}~\cite{Cockett2007}, 
 generalising the Van Kampen squares of~\cite{Lack2004}. 
 Indeed,
 examples of \textsc{vk}-colimits include
 coproducts in extensive categories and
 pushouts along monos in adhesive categories;
 the simplest example is a strict initial object. 

 The definition of \vk-colimits relies 
 only on elementary notions of category theory. 
 This feature, while attractive because of the implied simplicity, 
 obscures relationships 
 with other categorical concepts; the \emph{mathematical meaning}
 of \vk-colimits, so to speak.
 More abstract characterisations exist 
 for extensive and adhesive categories. 
 For instance, a category $\cat{C}$ is extensive if and only if the functor
 $+\colon\cat{C} \,{\downarrow}\, A \times \cat{C} \,{\downarrow}\,B\to \cat{C} \,{\downarrow}\, A+B$
 is an equivalence for any $A,B\in\cat{C}$~\cite{Lawvere1991,Carboni1993};
 adhesive categories can be characterised 
 in a similar manner~\cite{Lack2005}. 
 Our definition of \textsc{vk}-cocone is of the latter kind, 
 i.e.\ in terms of an equivalence of categories.
 We also provide an elementary characterisation in the spirit of Cockett and Guo.

 \smallskip
 This paper contains one central theorem: \textsc{vk}-cocones are
 those diagrams that are bicolimit diagrams when embedded in the associated
 bicategory of spans.  Bicolimits
 are the canonical notion of colimit in a bicategory.
 This characterises ``being Van Kampen'' as a universal
 property. We believe that this insight captures and explains the essence of the 
 various aforementioned well-behaved colimits studied in the literature.



 \paragraph{Structure of the paper.}
 In Section~\ref{sec:spans} we examine the category of spans and its applications.
 In Section~\ref{sec:prelude}
 we exhibit the relationship between coproducts in extensive categories and
 coproducts in their associated categories of spans, and explain why it is necessary
 to consider \emph{bicategories} of spans in order 
 to extend this relationship to arbitrary colimits.
 In Section~\ref{sec:preliminaries} we isolate the relevant class of bicategories 
 and recall the related notions. 
 In Section~\ref{sec:van-kampen-diagrams}
 we give a definition of \textsc{vk}-cocones together with an elementary characterisation
 and several examples. In Section~\ref{sec:vkinbicats} we recall the definition of
 bicolimits and prove several technical lemmas that allow us to pass between related
 concepts in a category and its associated bicategory of spans. The main theorem is
 proved in Section~\ref{sec:vanKampenSpans}. 


\section{Spans and generalised relations}
\label{sec:spans}



\noindent There are several concepts in category theory that generalise relations between sets.
The set theoretical concept of multirelation from $C$ to $D$
is a span of functions $C \gets[l] X \to[r] D$, which 
we will denote $\propto[X]{l}{r}\colon C \rightharpoonup D$, or simply \spa{l}{r}.
The set $X$ is sometimes referred to as the \emph{carrier}. 
Roughly, a pair of elements $(c,d) \in C \times D$ can be related in a number of ways;
concretely this is determined by 
the size of the preimage at $(c,d)$ of the function
$\langle l,r \rangle\colon X\rightarrow C\times D$, i.e.\  $\langle l,r \rangle^{-1}(c,d)$.
Two such spans
$\propto[X]{l}{r}$ and $\propto[X']{l'}{r'}$ would normally be considered equivalent if 
there is a bijection $\varphi\colon X\to X'$ that satisfies $l'\varphi=l$ and $r'\varphi=r$.
The existence of such a bijection yields an equivalence relation; quotienting by this
gives what we shall refer to as an \emph{abstract} span. Sometimes, for emphasis, we
shall refer to ordinary spans as \emph{concrete} spans. A multirelation is thus
an abstract span in the category of sets and functions.
Ordinary relations are captured by those abstract 
spans in which $\langle l,\,r\rangle$ is injective.

The concept of a span of morphisms makes sense in any category $\cat{C}$,
not only $\Set$; indeed we shall make use of the notation introduced in the previous
paragraph generally.
Extra structure is needed in order to compose spans, hence from this point onwards
we assume that $\CAT{C}$ has pullbacks.
Below we give an overview of the \emph{category} of $\CAT{C}$-spans, denoted by $\Spc{C}$.
Figure~\ref{fig:span-categoryFirstOccurrence:journal} accounts for its morphisms, 
identities and composition.
\begin{figure}[htb]
  \centering
         \begin{tikzpicture}[baseline={(0,0)}]
            \mnode[A]{C}{0,0}
            \mnode[A']{X
            }{1,0}
            \mnode{D}{2,0}
            \ardrawL{A'}{A}{above}{l}
            \ardrawL{A'}{D}{above}{r}
         \begin{scope}[>=left to,->,thick
           ]
           \draw (A.north) 
           ..controls +(.6,.6) and +(-.6,.6)..
           node[above] 
           {\ensuremath{\scriptstyle \spa{l}{r}}} 
           (D.north);
         \end{scope}
         \mnode[X']{X'
         }{1,-.6}
         \ardrawL{X'}{A}{below left}{l'}
         \ardrawL{X'}{D}{below right}{r'}
         \draw[draw=none] (X') -- (A') node[midway]{\rotatebox{-90}{$\cong$}}; 
          \end{tikzpicture}
          \qquad
  \begin{tikzpicture}[baseline={(0,0)}]
            \mnode{C}{0,0}
            \mnode[A']{C
            }{1,0}
            \mnode[B]{C}{2,0}
            \ardrawL{A'}{C}{above}{\text{id}}
            \ardrawL{A'}{B}{above}{\text{id}}
         \begin{scope}[>=left to,->,thick
           ]
           \draw (A.north) 
           ..controls +(.6,.6) and +(-.6,.6)..
           node[above] 
           {\ensuremath{\scriptstyle \spa{\id}{\id}}} 
           (B.north);
         \end{scope}
          \end{tikzpicture}
          \quad
  \begin{tikzpicture}[scale=1.7,baseline={(0,0)}]
          \tikzstyle{every node}+=[inner sep=1.5pt]
         \mnode[A]{C}{-.7,0}
         \mnode[A']{
         X}{0,0}
         \mnode[B]{D}{.5,-.5}
         \mnode[A'']{Z
         }{.5,.5}
         \mnode[B']{Y
         }{1,0}
         \mnode[C]{E}{1.7,0}
         \ardrawL{A'}{A}{above right}{l}
         \ardrawL{A'}{B}{above right,pos=.5}{r}
         \ardrawL{B'}{B}{above left,pos=.5}{l'}
         \ardrawL{B'}{C}{above left}{r'}
         \ardrawL{A''}{A'}{above left,pos=.5}{p}
         \ardrawL{A''}{B'}{above right,pos=.5}{q}
         \POC{[yshift=-4ex]A''.north}{-90} %
         \begin{scope}[>=left to,->,thick]
           \draw[
           ] (A) 
           ..controls +(.2,-.2) and +(-.7,0).. 
           node[pos=.4,below left] 
           {\ensuremath{\scriptstyle \spa{l}{r}}} 
           (B);
           \draw[
           ] (B)
           ..controls +(.7,0) and +(-.2,-.2).. 
           node[pos=.6,below right]
           {\ensuremath{\scriptstyle \spa{l'}{r'}}}
           (C);
           \draw[
           ] (A.north east)
           ..controls +(.9,.9) and +(-.9,.9).. 
           node[pos=.5,above]
           {\ensuremath{\scriptstyle \spa{lp}{r'q} }}
           (C.north west);
         \end{scope}
       \end{tikzpicture}
  \caption{Abstract spans, identities and composition via pullback}
  \label{fig:span-categoryFirstOccurrence:journal}
\end{figure}
More explicitly, $\Spc{C}$ has the same objects as $\CAT{C}$ and a 
morphism from $C$ to $D$ is an
equivalence class of spans $C \gets[l] X \to[r] D$.
 The identity on an object $C$ is a span of identities in $\CAT{C}$;
composition is obtained via pullback as illustrated in 
Figure~\ref{fig:span-categoryFirstOccurrence:journal} on the right.
It is not difficult to check that these definitions yield a category.
Moreover, when $\cat{C}=\Set$ this yields the expected composition
of multirelations. Moreover, $\CAT{C}$ has a standard embedding  $\Gr \colon \CAT{C} \to \Spc{C}$
into the span category.
This inclusion acts as the identity on objects and maps each morphism $f \colon A \to B$ in $\CAT{C}$
to its \emph{graph} $\Gr(f \colon A \to B) = \spa{\id_A}{f}$;
hence this embedding is also referred to as the \emph{graphing functor} 
\cite{Freyd1990}.

\begin{rem}\rm
Another way of generalising the concept of relation between sets to 
``relations'' between categories
is via the notion of \emph{profunctor} from $\CAT{C}$ to $\CAT{D}$:
it is an ordinary functor 
${\mathcal{F} \colon \CAT{C} \times  \CAT{D}^{op} \to \CAT{Set}}$ 
(here $\CAT{Set}$ could be replaced by another suitable monoidal category~$\mathcal{V}$).
Composition is via the left Kan extension along the Yoneda embedding~\cite{Borceux1994}.
The resulting structure is not a category but a bicategory.
Multirelations from $C$ to $D$ can be seen as profunctors where
$C$ and $D$ are considered as discrete categories. Formalising this observation yields
a biequivalence from $\Span[\Set]$ to the bicategory of profunctors between discrete 
categories~\cite{Lack2010}. This fact relies on a special property
of $\Set$, namely the equivalence of categories $\CAT{Set}\downarrow C \cong [C,\CAT{Set}]$,
and therefore does not generalise readily.
\end{rem}


Spans occur in very different contexts and 
often allow succinct characterisations of various concepts.
For example: Katis, Sabadini and Walters~\cite{Katis1997a} 
use spans to model systems with boundary (see also~\cite{Gadducci1999});
bisimulation has been captured as a span of open maps~\cite{Joyal1994}
as well as a span of coalgebra morphisms~\cite{Aczel1989};
an internal category is a monad in the bicategory of spans~\cite{B'enabou1967};
interaction categories~\cite{Abramsky1993} can be seen as examples of
process categories~\cite{Cockett1997}, which are certain
quotients of span bicategories; 
Mackey functors are coproduct preserving functors from the 
span-category~\cite{Lindner1976}.


\section{Colimits in the span category}
\label{sec:prelude}





\noindent We have seen that $\cat{C}$ embeds into $\Spc{C}$ and spans can be thought of as 
generalised relations. It is well-known (see, for example, \cite[1.911]{Freyd1990})
that in toposes colimits are preserved into the associated category of relations
via the standard embedding. A natural question then is what conditions 
of $\cat{C}$ colimits are necessary and sufficient for them to be
preserved by the embedding into $\Spc{C}$. We begin our investigation with the concrete case
of coproducts. It turns out that this problem is closely related with the notion of
extensive categories, where 
coproducts interact well with pullbacks.

The following, 
elementary definition of extensive categories makes this explicit~\cite{Carboni1993}.
\begin{figure}[htb]
  \centering
  $\displaystyle%
  \begin{tikzpicture}[scale=1.3,baseline=(current bounding box.west)]
    \mnode[AB]{A\,{+}\,B}{0,0}
    \mnode{A}{-1,0}
    \mnode{B}{1,0}
    \mnode{Z}{0,1}
    \mnode{X}{-1,1}
    \mnode{Y}{1,1}
    \ardrawL{A}{AB}{below}{\strt[]i_1}
    \ardrawL{B}{AB}{below}{\strt[]i_2}
    \ardrawL{X}{A}{left}{x}
    \ardrawL{Y}{B}{right}{y}
    \ardrawL{Z}{AB}{}{z}
    \ardrawL{X}{Z}{above}{\strt[]m}
    \ardrawL{Y}{Z}{above}{\strt[]n}
    \begin{scope}[white]
      \mnode[XY]{X\,{+}\,Y}{0,1.7}
      \ardrawL{X}{XY}{above left}{\strt[]j_1}
      \ardrawL{Y}{XY}{above right}{\strt[]j_2}
      \begin{scope}[draw=none]
        \draw[draw=none] 
        (XY) -- (Z) 
        node[midway]{
          \rotatebox{-90}{\ensuremath{\cong}}}
        ; 
      \end{scope}
    \end{scope}
  \end{tikzpicture}
  \quad\Rightarrow\quad
  \color{lightgray}\left(\quad
    \color{black}
      \begin{tikzpicture}[scale=1.3,baseline=(current bounding box.west)]
    \mnode[AB]{A\,{+}\,B}{0,0}
    \mnode{A}{-1,0}
    \mnode{B}{1,0}
    \mnode{Z}{0,1}
    \mnode{X}{-1,1}
    \mnode{Y}{1,1}
    \ardrawL{A}{AB}{below}{\strt[]i_1}
    \ardrawL{B}{AB}{below}{\strt[]i_2}
    \ardrawL{X}{A}{left}{x}
    \ardrawL{Y}{B}{right}{y}
    \ardrawL{Z}{AB}{}{z}
    \ardrawL{X}{Z}{above,pos=.7}{\strt[]m}
    \ardrawL{Y}{Z}{above,pos=.7}{\strt[]n}
    \begin{scope}[]
      \mnode[XY]{X\,{+}\,Y}{0,1.7}
      \ardrawL{X}{XY}{above left}{\strt[]j_1}
      \ardrawL{Y}{XY}{above right}{\strt[]j_2}
      \begin{scope}[draw=none]
        \draw[draw=none] (XY) -- (Z) node[midway]{\rotatebox{-90}{\ensuremath{\cong}}}; 
      \end{scope}
    \end{scope}
  \end{tikzpicture}
    \quad\Leftrightarrow  \quad
  \begin{tikzpicture}[scale=1.3,baseline=(current bounding box.west)]
    \mnode[AB]{A\,{+}\,B}{0,0}
    \mnode{A}{-1,0}
    \mnode{B}{1,0}
    \mnode{Z}{0,1}
    \mnode{X}{-1,1}
    \mnode{Y}{1,1}
    \ardrawL{A}{AB}{below}{\strt[]i_1}
    \ardrawL{B}{AB}{below}{\strt[]i_2}
    \ardrawL{X}{A}{left}{x}
    \ardrawL{Y}{B}{right}{y}
    \ardrawL{Z}{AB}{}{z}
    \ardrawL{X}{Z}{above}{\strt[]m}
    \ardrawL{Y}{Z}{above}{\strt[]n}
    \begin{scope}[white]
      \mnode[XY]{X\,{+}\,Y}{0,1.7}
      \ardrawL{X}{XY}{above left}{\strt[]j_1}
      \ardrawL{Y}{XY}{above right}{\strt[]j_2}
      \begin{scope}[draw=none]
        \draw[draw=none] (XY) -- (Z) node[midway]{\rotatebox{-90}{\ensuremath{\cong}}}; 
      \end{scope}
    \end{scope}
    \POC{X}{-45}
    \POC{Y}{-135}
  \end{tikzpicture}
    \color{lightgray}
    \quad\right)\color{black}
  $
  \caption{Extensivity condition}
  \label{fig:DefExtensivity:vks}
\end{figure}%
\begin{defi}
  A category 
  is extensive when
  \begin{enumerate}[(1)]
  \item it has finite coproducts; 
  \item it has pullbacks along coproduct injections;
  \item given a diagram on the left of Figure~\ref{fig:DefExtensivity:vks}
    with the bottom row a coproduct diagram,
    the top row is a coproduct diagram if only if the two squares are pullbacks.
  \end{enumerate}
\end{defi}\medskip

\noindent The relevant observation about coproducts in extensive categories is 
that the universal property of coproducts does not only
apply to morphisms of the category itself but actually extends to spans.
More precisely, 
given an ordinary coproduct $A \to[{i_1}] A + B \gets[{i_2}] B$
and a pair of morphisms $f \colon A \to C$ and $g \colon B \to C$ 
there exists unique mediating morphism $[f,g] \colon A + B \to C$;
given a coproduct in an extensive category and a pair of spans
$A \rightharpoonup[\spa{x}{h}] C$ and 
$B \rightharpoonup[\spa{y}{k}] C$
there exists a unique mediating span 
$A + B \rightharpoonup[\spa{x + y}{[h,k]}] C$.
This was already noticed by Lindner~\cite{Lindner1976}.

More can be said: it turns out that if a coproduct in $\cat{C}$
is also a coproduct in $\Spc{\cat{C}}$ then it satisfies the extensivity condition of
Figure~\ref{fig:DefExtensivity:vks}. 
Hence, the extensivity condition characterises the universal
property of coproducts in the ``larger universe'' of spans.
Summarizing, 
a coproduct in $\cat{C}$ satisfies the extensivity condition if and only if it
is a coproduct in $\Spc{\cat{C}}$.

\begin{prop}
  \label{prop:characterisationOfExtensivity:vks}
  Let $\CAT{C}$ be a category with coproducts and pullbacks.
  Then $\CAT{C}$ is extensive if and only 
  if the graphing functor~$\Gr \colon \CAT{C} \to \Spc{C}$ preserves coproducts.
\end{prop}
\begin{proof}
The fact that $\Gr$ preserves coproducts when $\cat{C}$ is extensive
was shown by Lindner~\cite[Lemma~3]{Lindner1976}.

For the converse, assume that $\Gr \colon \CAT{C} \to \Spc{C}$
preserves coproducts. First we shall show that coproducts in $\cat{C}$ are
stable under pullback. Let  $A \to[i_A] A +B \gets[i_B]B$
be a coproduct diagram and $z \colon Z \to A+B$.
Consider the diagram below.
  \[
  \xymatrix@=15pt{
  & {X+Y} \\ 
  {X} \ar[r]|{i} 
  \ar@/^1pc/[ur]^(.3){i_X}  
  \ar[d]_x & {Z} 
  \ar[u]_k 
  \ar[d]^	z & {Y} \ar[l]|j
  \ar@/_1pc/[ul]_(.3){i_Y} \ar[d]^y    
  \\
  {A} \ar[r]_-{i_A} & {A+B} & \ar[l]^(.3){i_B} {B}
  }
  \]
  First
  assume that
  $A \gets[x]X \to[i]Z$ and $Z \gets[j]Y \to[y]B$
  are pullbacks
  of $A \to[i_A] A+B \gets[z]Z$ and $Z \to[z]A+B\gets[i_B]B$,
  respectively. The existence of $k\colon Z\to X+Y$ with $ki=i_X$ and $kj=i_Y$	
  follows from the fact
  that $\langle\Gamma i_A,\,\Gamma i_B\rangle$ is a coproduct diagram in 
  $\Spc{\cat{C}}$.
  The universal property of $X+Y$ in $\cat{C}$ implies that 
  $k[i,j]=\id_{X+Y}$~${(*)}$. The universal property of $A+B$ in 
  $\Spc{C}$ implies that $\spa{z}{\id}=\spa{z}{[i,j]k}$ in
  $\Spc{C}$, which implies the existence of an isomorphism
  $\varphi\colon Z\to Z$ with $[i,j]k=\varphi$~${(**)}$. It now follows from
  ${(*)}$ that $k$ is split epi and from ${(**)}$ that it is mono; thus
  $k$ is an isomorphism.
  
  To verify the second part of the extensivity condition, consider the 
  boundary
  of the diagram above. We need to show that  
  $A \gets[x]X \to[i_X]X+Y$ and $X+Y \gets[i_Y]Y \to[y]B$ are pullbacks of
  $A \to[i_A] A+B \gets[x+y]X+Y$ and $X+Y \to[x+y]A+B\gets[i_B]B$.
  Now since $A+B$ is a coproduct in $\Spc{C}$ we get the existence of the interior
  part of the diagram, with the two squares pullbacks. By the argument in the
  previous paragraph $Z$ is the coproduct of $X$ and $Y$
  and so $k$ is an isomorphism.
\end{proof}


The main insight that can be
gained by inspection of this proof is a correspondence between
existence and uniqueness of mediating  spans on the one
hand and the two directions of the bi-implication of the Extensivity
Condition in Figure~\ref{fig:DefExtensivity:vks} on the other hand.
An analogous correspondence will recur later in the development of the main 
result in the bicategory of spans. 
The necessity of the bicategorical setting when considering arbitrary colimits
is the topic of the remainder of this section.

\subsection{The abstract span category is not enough}
\label{sec:why-span-category}

Proposition~\ref{prop:characterisationOfExtensivity:vks} could 
tempt one to try a generalisation to pushouts in the 
sense that a  pushout is ``well-behaved'' in $\cat{C}$
if and only if it is preserved by $\Gamma$ as a pushout in $\Spc{C}$.

A good candidate for a ``well-behaved'' pushout is given 
by the notion of Van Kampen square, which appeared as part of the
definition of adhesive categories~\cite{Lack2004}; indeed
this definition was partly motivated by the extensivity condition
of coproducts.

\begin{figure}[htb]
 \centering 
 \ensuremath{\displaystyle  
   \begin{tikzpicture}[baseline=(current bounding box.west),scale=1.6]
     \tikzstyle{every node}=[circle,inner sep=0pt]
     \mnode{B}{-1,0,0} ;
     \node (C) at (1,0,0) {$C$} ;
     \node (A) at (0,0,-1) {$A$} ;
     \node (D) at (0,0,1) {$D$} ;
     \node (B') at (-1,1,0) {$B'$} ;
     \node (C') at (1,1,0) {$C'$} ;
     \node (A') at (0,1,-1) {$A'$} ;
     \begin{scope}[lightgray,thick]
       \node (D') at (0,1,1) {$\boldsymbol{D'}$} ;
     \end{scope}
     \POCS {A'} {-1,0,+1} {0,-1,0}
     \POCS {A'} {+1,0,+1} {0,-1,0}
     \begin{scope}[thick]
       \ardrawL A B {below,very near start}{\strt f};
     \end{scope}
     \begin{scope}[thick]
       \ardrawL A C {below left,midway}{m\,};      
     \end{scope}
     \ardrawL {A'} {A} {right,near end}{\,a} ;
     \ardrawL{B'} {B} {left,near end}{b} ;   
     \ardrawL {C'} {C} {right,near end}{\,c};      
     \ardrawL {A'} {B'} {above,near end}{f'};
     \ardrawL {A'} {C'} {above right,midway}{\,m'};
     \begin{scope}[lightgray,thick]
       \ardrawL {B'} {D'} {below left, near end}{{\boldsymbol{n'}}^{\vphantom l}} ;
     \end{scope}
     \begin{scope}[thick]
       \ardrawL {B} {D} {below left,near end}{n\,};
     \end{scope}
     \begin{scope}[lightgray,thick]
       \ardrawLd {D'} {D} {left,pos=.68}{\boldsymbol{d}} ;
       \ardrawLd {C'} {D'} {below,near end}{\boldsymbol{g'}} ;   
     \end{scope}
     \begin{scope}[thick]
       \ardrawLd {C} {D} {below,near start}{g\strt};
     \end{scope}
   \end{tikzpicture}
   \Rightarrow
   \color{lightgray}
   \left(
     \color{black}
     \begin{tikzpicture}[baseline=(current bounding box.west),scale=1.6]
       \tikzstyle{every node}=[circle,inner sep=0pt]
       \begin{scope}[lightgray,thick]
         \mnode[B]{\boldsymbol{B}}{-1,0,0} ;
         \node (C) at (1,0,0) {$\boldsymbol C$} ;
         \node (A) at (0,0,-1) {$\boldsymbol A$} ;
         \node (D) at (0,0,1) {$\boldsymbol D$} ;
       \end{scope}
       \node (B') at (-1,1,0) {$B'$} ;
       \node (C') at (1,1,0) {$C'$} ;
       \node (A') at (0,1,-1) {$A'$} ;
       \node (D') at (0,1,1) {$D'$} ;
       \POCS[.17] {D'} {-1,0,-1} {+1,0,-1} 
       \begin{scope}[lightgray,thick]
         \ardrawL A B {below,very near start}{\strt \boldsymbol f};
         \ardrawL A C {below left,midway}{\boldsymbol m\,}; 
         \ardrawL {A'} {A} {right,near end}{\,\boldsymbol a} ;
         \ardrawL{B'} {B} {left,near end}{\boldsymbol b} ;   
         \ardrawL {C'} {C} {right,near end}{\,\boldsymbol c};      
       \end{scope}
       \ardrawL {A'} {B'} {above,near end}{f'};
       \ardrawL {A'} {C'} {above right,midway}{\,m'};
       \ardrawL {B'} {D'} {below left, near end}{{n'}^{\vphantom l}} ;   
       \begin{scope}[lightgray,thick]
         \ardrawL {B} {D} {below left,near end}{\boldsymbol n\,};   
         \ardrawLd {D'} {D} {left,pos=.68}{\boldsymbol d} ;
       \end{scope}
       \ardrawLd {C'} {D'} {below,near end}{g'} ;   
       \begin{scope}[lightgray]
         \ardrawLd {C} {D} {below,near start}{\boldsymbol g\strt}; 
       \end{scope}
     \end{tikzpicture}
     \Leftrightarrow
     \begin{tikzpicture}[baseline=(current bounding box.west),scale=1.6]
       \tikzstyle{every node}=[circle,inner sep=0pt]
       \mnode{B}{-1,0,0} ;
       \node (C) at (1,0,0) {$C$} ;
       \begin{scope}[lightgray]
         \node (A) at (0,0,-1) {$\boldsymbol A$} ;
       \end{scope}
       \node (D) at (0,0,1) {$D$} ;
       \node (B') at (-1,1,0) {$B'$} ;
       \node (C') at (1,1,0) {$C'$} ;
       \begin{scope}[lightgray]
         \node (A') at (0,1,-1) {$\boldsymbol {A'}$} ;
       \end{scope}
       \node (D') at (0,1,1) {$D'$} ;
       \POCS {C'} {-1,0,+1} {0,-1,0}
       \POCS {B'} {+1,0,+1} {0,-1,0}
       \begin{scope}[lightgray,thick]
         \ardrawL A B {below,very near start}{\strt \boldsymbol f}
         \ardrawL A C {below left,midway}{m\,} 
         \ardrawL {A'} {A} {right,near end}{\,a} 
       \end{scope}
       \ardrawL{B'} {B} {left,near end}{b}    
       \ardrawL {C'} {C} {right,near end}{\,c}      
       \begin{scope}[lightgray,thick]
         \ardrawL {A'} {B'} {above,near end}{\boldsymbol{f'}}
         \ardrawL {A'} {C'} {above right,midway}{\,\boldsymbol{m'}}
       \end{scope}
       \ardrawL {B'} {D'} {below left,near end}{{n'}^{\vphantom l}} 
       \ardrawL {B} {D} {below left,near end}{n\,}   
       \ardrawLd {D'} {D} {left,pos=.68}{d} 
       \ardrawLd {C'} {D'} {below,near end}{g'}    
       \ardrawLd {C} {D} {below,near start}{g\strt} 
     \end{tikzpicture}
     \color{lightgray}\right)\color{black}}
 \caption{Van Kampen square}
 \label{fig:vanKampen-once-more-span:diss}
\end{figure}
\begin{defi}[Van Kampen square]\rm
 \label{def:VK-span-later:diss}
 A commutative square $B \gets[f] A \to[m] C$, 
 $B \to[n]D \gets[g] C$  
 is \emph{Van Kampen} 
 when for each commutative cube
 as illustrated in Figure~\ref{fig:vanKampen-once-more-span:diss} on the left
 that has pullback squares as rear faces,
 its \emph{top~face is a pushout square} 
 if and only if
 its \emph{front faces are pullback squares}
 (see Figure~\ref{fig:vanKampen-once-more-span:diss}).
\end{defi}
A category is adhesive when it has pushouts along monomorphisms and these are Van Kampen
squares.

Differently from coproducts in extensive categories,
Van Kampen squares
do \emph{not} induce pushouts via inclusion
into the span category. Roughly, the reason for this is 
that for every pair of concrete $\CAT{C}$-spans $C \gets[l]X \to[r]D$
and $C \gets[l']X' \to[r']D$,
there may be several different isomorphisms $\varphi_1,\varphi_2,\ldots \colon X \to X'$
which witness that $\spa{l}{r}\colon C \rightharpoonup D$ and 
$\spa{l'}{r'}\colon A \rightharpoonup B$ are the same arrow in $\Spc{C}$.


\begin{figure}[htb]
  \centering
  \subfloat[\strut pushout in $\CAT{C}$]{%
    \label{subfig:Pushout}
    \begin{tikzpicture}[scale=1.5]
        \node[outer sep=0pt,inner sep=0pt] (D) at (0,-1) {%
    \begin{tikzpicture}[show background rectangle,
      background rectangle/.style={fill=black!30,rounded corners=0ex},scale=.5,outer sep=0pt,inner sep=.5pt]
    \node (a) at (0,0) {$\circ$};
    \node (b) at (1,0) {$\circ$};
    \node (ab) at (.5,-1) {$\star$};
    \draw[|->] (a) -- (ab);
    \draw[|->] (b) -- (ab);
    \end{tikzpicture}};

    \node[outer sep=0pt,inner sep=0pt] (C) at (-1,0,0) {
    \begin{tikzpicture}[show background rectangle,
      background rectangle/.style={fill=black!30,rounded corners=0ex},scale=.5,outer sep=0pt,inner sep=.5pt]
    \node (a) at (0,0) {$\circ$};
    \phantom{\node (b) at (1,0) {$\circ$};}
    \node (ab) at (.5,-1) {$\star$};
    \draw[|->] (a) -- (ab);
    \phantom{\draw[|->] (b) -- (ab);}
   \end{tikzpicture}};

    \node[outer sep=0pt,inner sep=0pt] (B) at (1,0,0) {
    \begin{tikzpicture}[show background rectangle,
      background rectangle/.style={fill=black!30,rounded corners=0ex},scale=.5,outer sep=0pt,inner sep=.5pt]
      \phantom{\node (a) at (0,0) {$\circ$};}
    \node (b) at (1,0) {$\circ$};
    \node (ab) at (.5,-1) {$\star$};
    \phantom{\draw[|->] (a) -- (ab);}
    \draw[|->] (b) -- (ab);
    \end{tikzpicture}};

    \node[outer sep=0pt,inner sep=0pt] (A) at (0,1) {
    \begin{tikzpicture}[show background rectangle,
      background rectangle/.style={fill=black!30,rounded corners=0ex},scale=.5,outer sep=0pt,inner sep=.5pt]
      \phantom{\node (a) at (0,0) {$\circ$};}
    \phantom{\node (b) at (1,0) {$\circ$};}
    \node (ab) at (.5,-1) {$\star$};
    \phantom{\draw[|->] (a) -- (ab);}
    \phantom{\draw[|->] (b) -- (ab);}
   \end{tikzpicture}};
 \ardrawL{A}{B}{above right}{m}
 \ardrawL{A}{C}{above left}{n}
 \ardrawL{C}{D}{below left}{m'} 
 \ardrawL{B}{D}{below right}{n'}
 \phantom{
    \node[outer sep=0pt,inner sep=0pt] at (0,-2) {
    \begin{tikzpicture}[show background rectangle,
      background rectangle/.style={fill=black!30,rounded corners=0ex},scale=.5,outer sep=0pt,inner sep=.5pt]
      \phantom{\node (a) at (0,0) {$\circ$};}
    \phantom{\node (b) at (1,0) {$\circ$};}
    \node (ab) at (.5,-1) {$\star$};
    \phantom{\draw[|->] (a) -- (ab);}
    \phantom{\draw[|->] (b) -- (ab);}
   \end{tikzpicture}};
}
  \end{tikzpicture}}
\quad\qquad\subfloat[\strut cocone in $\Spc{C}$]{%
    \label{subfig:spanCocone}
  \begin{tikzpicture}[scale=1.5]
    \begin{scope}[shift={(0,1)}]
          \node[outer sep=0pt,inner sep=0pt] (C) at (-1,0,0) {
    \begin{tikzpicture}[show background rectangle,
      background rectangle/.style={fill=black!30,rounded corners=0ex},scale=.5,outer sep=0pt,inner sep=.5pt]
    \node (a) at (0,0) {$\circ$};
    \phantom{\node (b) at (1,0) {$\circ$};}
    \node (ab) at (.5,-1) {$\star$};
    \draw[|->] (a) -- (ab);
    \phantom{\draw[|->] (b) -- (ab);}
   \end{tikzpicture}};

    \node[outer sep=0pt,inner sep=0pt] (B) at (1,0,0) {
    \begin{tikzpicture}[show background rectangle,
      background rectangle/.style={fill=black!30,rounded corners=0ex},scale=.5,outer sep=0pt,inner sep=.5pt]
      \phantom{\node (a) at (0,0) {$\circ$};}
    \node (b) at (1,0) {$\circ$};
    \node (ab) at (.5,-1) {$\star$};
    \phantom{\draw[|->] (a) -- (ab);}
    \draw[|->] (b) -- (ab);
    \end{tikzpicture}};

    \node[outer sep=0pt,inner sep=0pt] (A) at (0,1) {
    \begin{tikzpicture}[show background rectangle,
      background rectangle/.style={fill=black!30,rounded corners=0ex},scale=.5,outer sep=0pt,inner sep=.5pt]
      \phantom{\node (a) at (0,0) {$\circ$};}
    \phantom{\node (b) at (1,0) {$\circ$};}
    \node (ab) at (.5,-1) {$\star$};
    \phantom{\draw[|->] (a) -- (ab);}
    \phantom{\draw[|->] (b) -- (ab);}
   \end{tikzpicture}};
 \pardrawL{A}{B}{above right}{\Gr(m)}
 \PardrawL{A}{C}{above left}{\Gr(n)}
    \end{scope}
    \node[outer sep=0pt,inner sep=0pt] (C') at (-1,0) {
    \begin{tikzpicture}[show background rectangle,
      background rectangle/.style={fill=black!30,rounded corners=0ex},scale=.5,outer sep=0pt,inner sep=.5pt]
    \node (a) at (0,0) {$\circ$};
    \phantom{\node (b) at (1,0) {$\circ$};}
    \node (a') at (0,-1) {$\star$};
    \node (b') at (1,-1) {$\star$};
    \draw[|->] (a) -- (a');
    \phantom{\draw[|->] (b) -- (b');}
    \end{tikzpicture}};

    \node[outer sep=0pt,inner sep=0pt] (B') at (1,0) {
    \begin{tikzpicture}[show background rectangle,
      background rectangle/.style={fill=black!30,rounded corners=0ex},scale=.5,outer sep=0pt,inner sep=.5pt]
      \phantom{\node (a) at (0,0) {$\circ$};}
    \node (b) at (1,0) {$\circ$};
    \node (a') at (0,-1) {$\star$};
    \node (b') at (1,-1) {$\star$};
    \phantom{\draw[|->] (a) -- (a');}
    \draw[|->] (b) -- (b');
    \end{tikzpicture}};

    \node[outer sep=0pt,inner sep=0pt] (T) at (0,-1) {
    \begin{tikzpicture}[show background rectangle,
      background rectangle/.style={fill=black!30,rounded corners=0ex},scale=.5,outer sep=0pt,inner sep=.5pt]
      \phantom{\node at (0,0) {$\circ$};}
      \node (t) at (.5,0) {$\circ$};
      \phantom{\node  at (1,0) {$\circ$};}
      \node (t') at (.5,-1) {$\star$};
      \draw[|->] (t) -- (t');
   \end{tikzpicture}};
 \begin{scope}[very thick,lightgray]
 \ardrawL{B'}{B}{right,pos=.3}{!}
 \ardrawL{C'}{C}{left,pos=.3}{!}
 \ardrawL{B'}{T}{below right,pos=.3}{!}
 \ardrawL{C'}{T}{below left,pos=.3}{!}
\end{scope}
\pardrawL{[xshift=1ex]B.south west}{T}{above left}{s}
 \PardrawL{[xshift=-1ex]C.south east}{T}{above right}{t}
  \end{tikzpicture}}
\quad\qquad\subfloat[\strut diagonal]{%
    \label{subfig:diagonal}
    \quad \begin{tikzpicture}[scale=1.5]
    \begin{scope}[shift={(0,1)}]

    \node[outer sep=0pt,inner sep=0pt] (A) at (0,1) {
    \begin{tikzpicture}[show background rectangle,
      background rectangle/.style={fill=black!30,rounded corners=0ex},scale=.5,outer sep=0pt,inner sep=.5pt]
      \phantom{\node (a) at (0,0) {$\circ$};}
    \phantom{\node (b) at (1,0) {$\circ$};}
    \node (ab) at (.5,-1) {$\star$};
    \phantom{\draw[|->] (a) -- (ab);}
    \phantom{\draw[|->] (b) -- (ab);}
   \end{tikzpicture}};

    \node[outer sep=0pt,inner sep=0pt] (U) at (0,-.5) {
    \begin{tikzpicture}[show background rectangle,
      background rectangle/.style={fill=black!30,rounded corners=0ex},scale=.5,outer sep=0pt,inner sep=.5pt]
      \phantom{\node (a) at (0,0) {$\circ$};}
    \phantom{\node (b) at (1,0) {$\circ$};}
    \node (a') at (0,-1) {$\star$};
    \node (b') at (1,-1) {$\star$};
    \phantom{\draw[|->] (a) -- (a');}
    \phantom{\draw[|->] (b) -- (b');}
   \end{tikzpicture}};
    \end{scope}

    \node[outer sep=0pt,inner sep=0pt] (T) at (0,-1) {
    \begin{tikzpicture}[show background rectangle,
      background rectangle/.style={fill=black!30,rounded corners=0ex},scale=.5,outer sep=0pt,inner sep=.5pt]
      \phantom{\node at (0,0) {$\circ$};}
      \node (t) at (.5,0) {$\circ$};
      \phantom{\node  at (1,0) {$\circ$};}
      \node (t') at (.5,-1) {$\star$};
      \draw[|->] (t) -- (t');
   \end{tikzpicture}};
 \begin{scope}[very thick,lightgray]
 \ardrawL{U}{A}{right,pos=.3}{!}
 \ardrawL{U}{T}{right,pos=.3}{!}
 \end{scope}
 \begin{scope}[bend left,looseness=.5]
 \pardrawL[to]{A.east}{T.east}{right}{d}
 \end{scope}
  \end{tikzpicture}\qquad} 
  \caption{Van Kampen square in $\CAT{C} = [ {\cdot}{\to}{\cdot}, \CAT{Set}]$ and a cocone of spans}
  \label{fig:ExampleOne}
\end{figure}
Our counterexample is a pushout in the category
$\CAT{Set}$-arrows,
i.e.\ the functor category $[ {\cdot}{\to}{\cdot}, \CAT{Set}]$,
which is adhesive \cite{Lack2005}.
The objects of this category are functions.  They will be  depicted as gray boxes in which all input-output pairs are connected by arrows of the form `$\mapsto$' (see Figure~\ref{fig:ExampleOne});
each element of the domain is rendered as a small circle `$\circ$', 
while each element of the codomain is represented by a star `$\star$'.


Now Figure~\ref{fig:ExampleOne}{\textsc{\subref{subfig:Pushout}} is a pushout of a pair of monomorphisms and hence
a Van Kampen square.
Consider the cocone for the ``same'' span of morphisms
in $\Spc{C}$ described in Figure~\ref{fig:ExampleOne}\textsc{\subref{subfig:spanCocone}}; in the latter 
figure the gray arrows belong to  $\CAT{C}$
while the arrows $s$ and $t$ belong to  $\Spc{C}$. 
Figure~\ref{fig:ExampleOne}\textsc{\subref{subfig:diagonal}} shows the diagonal span $d$ of this square.
Observe that any concrete representative of the diagonal span $d$ actually has a 
non-trivial symmetry group which consists of
 one non-identity automorphism as well as the identity.

\begin{figure}[htb]
  \centering
  \begin{tikzpicture}[scale=2.1,yscale=1.2]
    \node[outer sep=0pt,inner sep=0pt] (T) at (0,2.3,1) {
    \begin{tikzpicture}[show background rectangle,
      background rectangle/.style={fill=black!30,rounded corners=0ex},scale=.5,outer sep=0pt,inner sep=.5pt]
      \phantom{\node at (0,0) {$\circ$};}
      \node (t) at (.5,0) {$\circ$};
      \phantom{\node  at (1,0) {$\circ$};}
      \node (t') at (.5,-1) {$\star$};
      \draw[|->] (t) -- (t');
   \end{tikzpicture}};

    \node[outer sep=0pt,inner sep=0pt] (D) at (0,0,1) {
    \begin{tikzpicture}[show background rectangle,
      background rectangle/.style={fill=black!30,rounded corners=0ex},scale=.5,outer sep=0pt,inner sep=.5pt]
    \node (a) at (0,0) {$\circ$};
    \node (b) at (1,0) {$\circ$};
    \node (ab) at (.5,-1) {$\star$};
    \draw[|->] (a) -- (ab);
    \draw[|->] (b) -- (ab);
    \end{tikzpicture}};

    \node[outer sep=0pt,inner sep=0pt] (C) at (-1,0,0) {
    \begin{tikzpicture}[show background rectangle,
      background rectangle/.style={fill=black!30,rounded corners=0ex},scale=.5,outer sep=0pt,inner sep=.5pt]
    \node (a) at (0,0) {$\circ$};
    \phantom{\node (b) at (1,0) {$\circ$};}
    \node (ab) at (.5,-1) {$\star$};
    \draw[|->] (a) -- (ab);
    \phantom{\draw[|->] (b) -- (ab);}
   \end{tikzpicture}};

    \node[outer sep=0pt,inner sep=0pt] (B) at (1,0,0) {
    \begin{tikzpicture}[show background rectangle,
      background rectangle/.style={fill=black!30,rounded corners=0ex},scale=.5,outer sep=0pt,inner sep=.5pt]
      \phantom{\node (a) at (0,0) {$\circ$};}
    \node (b) at (1,0) {$\circ$};
    \node (ab) at (.5,-1) {$\star$};
    \phantom{\draw[|->] (a) -- (ab);}
    \draw[|->] (b) -- (ab);
    \end{tikzpicture}};

    \node[outer sep=0pt,inner sep=0pt] (A) at (0,0,-1) {
    \begin{tikzpicture}[show background rectangle,
      background rectangle/.style={fill=black!30,rounded corners=0ex},scale=.5,outer sep=0pt,inner sep=.5pt]
      \phantom{\node (a) at (0,0) {$\circ$};}
    \phantom{\node (b) at (1,0) {$\circ$};}
    \node (ab) at (.5,-1) {$\star$};
    \phantom{\draw[|->] (a) -- (ab);}
    \phantom{\draw[|->] (b) -- (ab);}
   \end{tikzpicture}};

 \begin{scope}[shift={(0,1)}]
    \node[outer sep=0pt,inner sep=0pt] (D') at (0,0,1) {
    \begin{tikzpicture}[show background rectangle,
      background rectangle/.style={fill=black!30,rounded corners=0ex},scale=.5,outer sep=0pt,inner sep=.5pt]
    \node (a) at (0,0) {$\circ$};
    \node (b) at (1,0) {$\circ$};
    \node (a') at (0,-1) {$\star$};
    \node (b') at (1,-1) {$\star$};
    \draw[|->] (a) -- (a');
    \draw[|->] (b) -- (b');
    \end{tikzpicture}};

    \node[outer sep=0pt,inner sep=0pt] (C') at (-1,0,0) {
    \begin{tikzpicture}[show background rectangle,
      background rectangle/.style={fill=black!30,rounded corners=0ex},scale=.5,outer sep=0pt,inner sep=.5pt]
    \node (a) at (0,0) {$\circ$};
    \phantom{\node (b) at (1,0) {$\circ$};}
    \node (a') at (0,-1) {$\star$};
    \node (b') at (1,-1) {$\star$};
    \draw[|->] (a) -- (a');
    \phantom{\draw[|->] (b) -- (b');}
    \end{tikzpicture}};

    \node[outer sep=0pt,inner sep=0pt] (B') at (1,0,0) {
    \begin{tikzpicture}[show background rectangle,
      background rectangle/.style={fill=black!30,rounded corners=0ex},scale=.5,outer sep=0pt,inner sep=.5pt]
      \phantom{\node (a) at (0,0) {$\circ$};}
    \node (b) at (1,0) {$\circ$};
    \node (a') at (0,-1) {$\star$};
    \node (b') at (1,-1) {$\star$};
    \phantom{\draw[|->] (a) -- (a');}
    \draw[|->] (b) -- (b');
    \end{tikzpicture}};

    \node[outer sep=0pt,inner sep=0pt] (A') at (0,0,-1) {
    \begin{tikzpicture}[show background rectangle,
      background rectangle/.style={fill=black!30,rounded corners=0ex},scale=.5,outer sep=0pt,inner sep=.5pt]
      \phantom{\node (a) at (0,0) {$\circ$};}
      \phantom{\node (b) at (1,0) {$\circ$};}
    \node (a') at (0,-1) {$\star$};
    \node (b') at (1,-1) {$\star$};
    \phantom{\draw[|->] (a) -- (a');}
    \phantom{\draw[|->] (b) -- (b');}
    \end{tikzpicture}};
 \end{scope}
 \ardraw{A}{B}
 \ardraw{A}{C}
 \ardraw{B}{D}
 \ardraw{C}{D}
 \ardraw{A'}{A}
 \ardraw{B'}{B}
 \ardraw{C'}{C}
 \begin{scope}[dashed]
 \ardrawd{D'}{D}
\end{scope}
\ardraw{A'}{B'}
 \ardraw{A'}{C'}
 \ardrawd{B'}{D'}
 \ardrawd{C'}{D'}
 \begin{scope}[dashed]
   \ardrawd{D'}{T}
 \end{scope}
 \ardraw[..controls +(0,.5) and +(.5,0,-.5)..]{B'}{T.east}
 \ardraw[..controls +(0,.5) and +(-.5,0,-.5)..]{C'}{T.west}
 \draw[draw=black!70,very thick] ([shift={(-.5ex,.5ex)}]T.north west) rectangle ([shift={(.5ex,-.5ex)}]D.south east);
 \draw[thick] ([shift={(.1ex,-.1ex)}]B'.north west) rectangle ([shift={(-.1ex,.1ex)}]B'.south east);
  \end{tikzpicture}
  \qquad  \qquad
  \begin{tikzpicture}[scale=2.1,yscale=1.2]
        \node[outer sep=0pt,inner sep=0pt] (T) at (0,2.3,1) {
    \begin{tikzpicture}[show background rectangle,
      background rectangle/.style={fill=black!30,rounded corners=0ex},scale=.5,outer sep=0pt,inner sep=.5pt]
      \phantom{\node at (0,0) {$\circ$};}
      \node (t) at (.5,0) {$\circ$};
      \phantom{\node  at (1,0) {$\circ$};}
      \node (t') at (.5,-1) {$\star$};
      \draw[|->] (t) -- (t');
   \end{tikzpicture}};

    \node[outer sep=0pt,inner sep=0pt] (D) at (0,0,1) {
    \begin{tikzpicture}[show background rectangle,
      background rectangle/.style={fill=black!30,rounded corners=0ex},scale=.5,outer sep=0pt,inner sep=.5pt]
    \node (a) at (0,0) {$\circ$};
    \node (b) at (1,0) {$\circ$};
    \node (ab) at (.5,-1) {$\star$};
    \draw[|->] (a) -- (ab);
    \draw[|->] (b) -- (ab);
    \end{tikzpicture}};

    \node[outer sep=0pt,inner sep=0pt] (C) at (-1,0,0) {
    \begin{tikzpicture}[show background rectangle,
      background rectangle/.style={fill=black!30,rounded corners=0ex},scale=.5,outer sep=0pt,inner sep=.5pt]
    \node (a) at (0,0) {$\circ$};
    \phantom{\node (b) at (1,0) {$\circ$};}
    \node (ab) at (.5,-1) {$\star$};
    \draw[|->] (a) -- (ab);
    \phantom{\draw[|->] (b) -- (ab);}
   \end{tikzpicture}};

    \node[outer sep=0pt,inner sep=0pt] (B) at (1,0,0) {
    \begin{tikzpicture}[show background rectangle,
      background rectangle/.style={fill=black!30,rounded corners=0ex},scale=.5,outer sep=0pt,inner sep=.5pt]
      \phantom{\node (a) at (0,0) {$\circ$};}
    \node (b) at (1,0) {$\circ$};
    \node (ab) at (.5,-1) {$\star$};
    \phantom{\draw[|->] (a) -- (ab);}
    \draw[|->] (b) -- (ab);
    \end{tikzpicture}};

    \node[outer sep=0pt,inner sep=0pt] (A) at (0,0,-1) {
    \begin{tikzpicture}[show background rectangle,
      background rectangle/.style={fill=black!30,rounded corners=0ex},scale=.5,outer sep=0pt,inner sep=.5pt]
      \phantom{\node (a) at (0,0) {$\circ$};}
    \phantom{\node (b) at (1,0) {$\circ$};}
    \node (ab) at (.5,-1) {$\star$};
    \phantom{\draw[|->] (a) -- (ab);}
    \phantom{\draw[|->] (b) -- (ab);}
   \end{tikzpicture}};
 \begin{scope}[shift={(0,1)}]
    \node[outer sep=0pt,inner sep=0pt] (D') at (0,0,1) {
    \begin{tikzpicture}[show background rectangle,
      background rectangle/.style={fill=black!30,rounded corners=0ex},scale=.5,outer sep=0pt,inner sep=.5pt]
    \node (a) at (0,0) {$\circ$};
    \node (b) at (1,0) {$\circ$};
    \node (a') at (0,-1) {$\star$};
    \node (b') at (1,-1) {$\star$};
    \draw[|->] (a) -- (a');
    \draw[|->] (b) -- (a');
    \end{tikzpicture}};

    \node[outer sep=0pt,inner sep=0pt] (C') at (-1,0,0) {
    \begin{tikzpicture}[show background rectangle,
      background rectangle/.style={fill=black!30,rounded corners=0ex},scale=.5,outer sep=0pt,inner sep=.5pt]
    \node (a) at (0,0) {$\circ$};
    \phantom{\node (b) at (1,0) {$\circ$};}
    \node (a') at (0,-1) {$\star$};
    \node (b') at (1,-1) {$\star$};
    \draw[|->] (a) -- (a');
    \phantom{\draw[|->] (b) -- (b');}
    \end{tikzpicture}};

    \node[outer sep=0pt,inner sep=0pt] (B') at (1,0,0) {
    \begin{tikzpicture}[show background rectangle,
      background rectangle/.style={fill=black!30,rounded corners=0ex},scale=.5,outer sep=0pt,inner sep=.5pt]
      \phantom{\node (a) at (0,0) {$\circ$};}
    \node (b) at (1,0) {$\circ$};
    \node (a') at (0,-1) {$\star$};
    \node (b') at (1,-1) {$\star$};
    \phantom{\draw[|->] (a) -- (a');}
    \draw[|->] (b) -- (a');
    \end{tikzpicture}};

    \node[outer sep=0pt,inner sep=0pt] (A') at (0,0,-1) {
    \begin{tikzpicture}[show background rectangle,
      background rectangle/.style={fill=black!30,rounded corners=0ex},scale=.5,outer sep=0pt,inner sep=.5pt]
      \phantom{\node (a) at (0,0) {$\circ$};}
      \phantom{\node (b) at (1,0) {$\circ$};}
    \node (a') at (0,-1) {$\star$};
    \node (b') at (1,-1) {$\star$};
    \phantom{\draw[|->] (a) -- (a');}
    \phantom{\draw[|->] (b) -- (b');}
    \end{tikzpicture}};
 \end{scope}
 \ardraw{A}{B}
 \ardraw{A}{C}
 \ardraw{B}{D}
 \ardraw{C}{D}
 \ardraw{A'}{A}
 \ardraw{B'}{B}
 \ardraw{C'}{C}
 \begin{scope}[dashed]
   \ardrawd{D'}{D}   
 \end{scope}

 \ardraw{A'}{B'}
 \ardraw{A'}{C'}
 \ardrawd{B'}{D'}
 \ardrawd{C'}{D'}
 \begin{scope}[dashed]
   \ardrawd{D'}{T}
 \end{scope}
 \ardraw[..controls +(0,.5) and +(.5,0,-.5)..]{B'}{T.east}
 \ardraw[..controls +(0,.5) and +(-.5,0,-.5)..]{C'}{T.west}
 \draw[draw=black!70,very thick] ([shift={(-.5ex,.5ex)}]T.north west) rectangle ([shift={(.5ex,-.5ex)}]D.south east);
 \draw[thick] ([shift={(.1ex,-.1ex)}]B'.north west) rectangle ([shift={(-.1ex,.1ex)}]B'.south east);
\end{tikzpicture}

  \caption{Two different mediating spans to the same cocone of spans}
  \label{fig:Counterexample}
\end{figure}

As illustrated in Figure~\ref{fig:Counterexample}, 
each of the mediating spans can be constructed in the
category~$\CAT{C}$ 
by means of the Van Kampen square property.
The cocone from Figure~\ref{fig:ExampleOne}\textsc{\subref{subfig:spanCocone}} 
consisting of $s$ and $t$ is now pointing to the top. 
Commutativity of the cocone in $\Spc{C}$ gives a pair of $\CAT{C}$-pullbacks in the back of each of the diagrams in Figure~\ref{fig:Counterexample}. 
 Depending on choices of pullbacks and  witnesses for the abstract equality of the composite spans,
we obtain two different diagrams over which we can take a pushout.
(In Figure~\ref{fig:Counterexample},
we tried to express this fact by ``switching'' of the two stars in the carrier of right span of the cocone.)
These two diagrams on top of a cube
yield two different mediating arrows in the category of spans by 
taking the pushout in $\CAT{C}$.
In the end, we obtain two different mediating spans to the same cocone in $\Spc{C}$, 
and thus $\cdot \rightharpoonup [\Gr(n')] \cdot \leftharpoonup[\Gr(m')]\cdot$ 
cannot be a pushout of $\cdot \leftharpoonup[\Gr(m)] \cdot \rightharpoonup[\Gr(n)] \cdot $.



\begin{rem}\rm
To solve the problem of non-trivial symmetry groups of spans one could try to restrict
to partial map spans, i.e.\ those spans with a monomorphisms from the carrier into the domain.
This however would yield a (properly) weaker notion,
which one could call partial Van Kampen square~\cite{Heindel2010a}.
%
%
In fact, in the category of sets and, more generally, any elementary topos, all colimits are
partial Van Kampen colimits while there are examples of pushouts in
Set, which are not Van Kampen~\cite{Lack2005}. 
\end{rem}

 The only canonical way to ``tame'' the
non-trivial symmetry groups of spans 
is to keep track of the involved isomorphisms,
i.e.\ we have  to work in a bicategory $\Span[C]$ of (concrete) spans in $\cat{C}$.

It turns out that by moving to the setting of bicategories (recalled in the proceeding
section) we obtain
not only that Van Kampen squares characterise those pushouts that are preserved
by $\Gamma$ but that a general Van Kampen condition
(introduced in Section~\ref{sec:van-kampen-diagrams}) characterises those colimits that
satisfy the bicategorical universal property of colimits.


\section{Bicategories}
\label{sec:preliminaries}


 \noindent Here we introduce background on bicategories~\cite{B'enabou1967,Borceux1994,Lack2010}
 and some notational conventions.
 Our focus is the bicategory of spans over a category $\CAT{C}$ with 
 a choice of pullbacks (see Example~\ref{exm:spans}).
 This allows us to avoid unnecessary book-keeping by considering
 only those bicategories that satisfy the identity axioms strictly.
 It is only a cosmetic choice,
 the development can be easily adapted to the 
 standard setting of bicategories.

 \begin{defi}[Strictly unitary bicategories]\rm
   \label{defn:subicat}
   A \emph{strictly unitary \textsc{(su)} bicategory} $\bicat{B}$ consists of:
   \begin{enumerate}[$\bullet$]
   \item a collection $\Ob \bicat{B}$ of \emph{objects};
   \item for $A,B\in \Ob \bicat{B}$ a category $\bicat{B}(A,B)$, the objects
     and arrows of which are called, respectively,
     the \emph{arrows} and the \emph{2-cells} of $\bicat{B}$.
     Composition is denoted by $\circledcirc$ and referred
     to as \emph{vertical composition}.
     Given $(f \colon A \to B) \in \bicat{B}(A,B)$, 
     its identity 2-cell will be denoted $\iota_f \colon f \to f$.
     Each $\bicat{B}(A,A)$ contains 
     a special object $\id_A \colon A \to A$, 
     called the \emph{identity arrow};
   \item for $A,B,C\in\Ob\bicat{B}$, 
     a functor
     $c_{A,B,C} \colon \bicat{B}(A,B) \times \bicat{B}(B,C)  \to \bicat{B}(A,C)$
     called \emph{horizontal composition}. 
     On objects, $c_{A,B,C}\langle f,g\rangle$ is written $g \circ f$, 
     while on arrows
     $c_{A,B,C}\langle\gamma,\delta\rangle$ it is
     $\delta \ast \gamma$.
     For any $f \colon A \to B$ we have 
     $\id_B\circ f = f =f\circ \id_A$;
   \item for $A,B,C,D\in\Ob\bicat{B}$, 
     arrows $f \colon A \to B$, 
     $g \colon B \to C$ and 
     $h \colon C \to D$ an \emph{associativity natural isomorphism}
     $\alpha_{A,B,C,D}(f,g,h) \colon 
     h\circ (g\circ f) \to (h\circ g)\circ f$.
     It satisfies the coherence axioms:
     for any composable $f,g,h,k$, we have
     $\alpha_{f,\id,g} = \iota_{g \circ f}$ 
     and also that the following {2-cells} are equal:
     \begin{displaymath}
       \qquad\begin{tikzpicture}[scale=1.5,yscale=1.3,baseline=(A.base)]
         \mnode{E}{0,0}
         \mnode{C}{2,.5}
         \mnode{B}{3,.5}
         \mnode{A}{4,0}
         \begin{scope}[bend right=15]
           \ardrawL[to]{C}{E.north}{above left}{k \circ h}
           \ardrawL{B}{C}{above}{g}
           \ardrawL[to]{A.north}{B}{above right}{f}
         \end{scope}
         \mnode[D']{D}{1,0}
         \mnode[C']{C}{2,0}
         \ardrawL{D'}{E}{above}{k}
         \ardrawL{C'}{D'}{above}{h}
         \ardrawL{A}{C'}{above}{g \circ f}
         \mnode[D'']{D}{1,-.5}
         \mnode[C'']{}{2,-.5}
         \ardrawL[{to[bend left=15]}]{D''}{E.south}{below left}{k}
         \ardrawL[{to[bend left=15]}]{A.south}{D''}{below right}%
         {h \circ (g \circ f)}
         \artwoL{C'.north}{C.south}{right}{\alpha}
         \artwoL{C''.north}{C'.south}{right}{\alpha}
       \end{tikzpicture}
       =
       \begin{tikzpicture}[scale=1.5,yscale=1.3,baseline=(A.base)]
         \mnode{E}{0,0}
         \mnode[D]{}{1,.5}
         \mnode{C}{2,.5}
         \mnode{B}{3,.5}
         \mnode[A]{A\,\,.\!\!\!}{4,0}
         \begin{scope}[bend right=15]
           \ardrawL[to]{C}{E.north east}{above left}{k \circ h}
           \ardrawL{B}{C}{above}{g}
           \ardrawL[to]{A.north}{B}{above right}{f}
         \end{scope}
         \mnode[D']{D}{1,0}
         \mnode[C']{}{2,0}
         \ardrawL{D'}{E}{above}{k}
         \ardrawL[{to[bend left=15]}]{B}{D'}{above left,pos=.4}{h \circ g}
         \mnode[D'']{D}{1,-.5}
         \ardrawL[{to[bend right=5]}]{A}{D''}{above,pos=.3}%
         {(h \circ g) \circ f_{\vphantom{|}}}
         \draw[draw=none](A) to[bend right=5](D'') coordinate[pos=.5](mid);

         \mnode[C'']{}{2,-.5}
         \ardrawL[{to[bend left=15]}]{D''}{E}{below left}{k}
         \ardrawL[{to[bend left=15]}]{A.south}{D''}{below right}%
         {h \circ (g \circ f)}
         \draw[draw=none](A.south) to [bend left=15] %
         (D'') coordinate[pos=.5](mid');
         \begin{scope}[shorten <=.5ex,shorten >=.5ex]
           \artwoL{D'.north east}{D.south east}{right}{\alpha}
         \end{scope}

         \artwoL{D''.north east}{[xshift=-.5ex]D'.south east}{right}{\alpha}
         \begin{scope}[shorten <=.5ex,shorten >=.5ex]
           \artwoL{mid'}{mid}{right}{\alpha}   
         \end{scope}

       \end{tikzpicture}
     \end{displaymath}

   \end{enumerate}
 \end{defi}

 \begin{exa}\rm
 Any (ordinary) category $\cat{C}$ is
 a (\textsc{su}-)bicategory with trivial 2-cells. 
 \end{exa}


 As we have seen, composition in the category \Spc{C} is obtained via pullback.
 Because of the universal property of pullbacks and the fact that the arrows
 of \Spc{C} are abstract spans, composition defined in this fashion is well-defined.
 Instead, in order to compose concrete spans we shall need to assume some \emph{choice}
 of pullback in $\cat{C}$; this means that for any cospan $X \to[f] Z \gets[g]Y$
 there is a chosen object $X \times_Z Y$ and span $X \gets[\pi_1] X \times_Z Y \to[\pi_2] Y$
 that together with $f$ and $g$ forms a pullback square.
 Moreover we assume that the choice preserves identities:
 if $f$ is $\id_X$ then $X\times_Z Y=Y$ and 
 $\pi_1=\id_Y$, and analogously for $g$. This is a completely harmless assumption since the
 identity of the chosen pullback diagram for any cospan is insignificant: any two 
 choices are equivalent. 

 \begin{exa}[Span bicategory~\cite{B'enabou1967}]\rm
   \label{exm:spans}   
 %
   Assume that $\cat{C}$  
   has a choice of pullbacks that preserves identities. 
   $\Span[C]$ has:
   \begin{enumerate}[$\bullet$]
   \item as objects, the objects of $\cat{C}$,
     i.e.\ $\Ob\Span[C] = \Ob{\CAT{C}}$;
   \item  as arrows from $C$ to $D$,
     the $\CAT{C}$-spans $C\gets[l]X \to[r] D$.
   The composition with another 
   span $D \gets[l'] Y \to[r']E$ is obtained via the chosen
   pullback as illustrated below;
   however this composition is only associative up to canonical isomorphism.
   The identity on an object~$C$ is the span $C \gets[\id] C \to[\id]C$.
   \[
  \begin{tikzpicture}[scale=1,baseline={(0,.9)}]
     \mnode [A]{C}{-.41,0}
     \mnode{X} {1,0}
     \mnode [B]{D} {2,-1}
     \mnode[U]{X\times_D Y} {2,1}
     \POC[.16]{U}{-90}
     \mnode [Z]{Y} {3,0}
     \mnode [C] {E} {4.41,0}
     \ardrawL X A {above} l
     \begin{scope}[shorten <=-1pt,shorten >=-1pt]
       \ardrawL X B {below left,pos=.2} r
       \ardrawL Z B {below right,pos=.2} {l'}
     \end{scope}
     \ardrawL Z C {above} {r'}
     \begin{scope}[shorten <=-1pt,shorten >=-1pt]
       \ardrawL U X {above left} {\pi_1}
       \ardrawL U Z {above right} {\pi_2}
     \end{scope}
     \PardrawL[..controls+(0,-1) and +(-.5,0)..]%
     {A}{B}{below,pos=.5}{\spa{l}{r}}
     \PardrawL[..controls+(+.5,0) and +(0,-1) ..]%
     {B}{C}{below,pos=.5}{\spa{l'}{r'}}
     \pardrawL[..controls+(1,1.8) and +(-1,1.8)..]%
     {A}{C}{above,pos=.5}{\spa{l\pi_1}{r'\pi_2}}
   \end{tikzpicture}%
   \]
  \item its 2-cells 
     $\xi \colon \spa{l}{r} \to \spa{l'}{r'}$ 
     are $\CAT{C}$-arrows
     $\xi \colon X \to X'$ between the respective carriers 
     such that $l' \circ \xi = l$ and 
     $r' \circ\xi = r$.
   \end{enumerate}
 \end{exa}

\noindent  For our purposes it suffices to consider \emph{strict homomorphisms} between
 \textsc{su}-bicategories.
 \begin{defi}[Strict homomorphisms~\cite{B'enabou1967}]%
   \rm\label{defn:strictHom}
   Let $\bicat{A}$ and $\bicat{B}$ be \textsc{su}-bicategories. A 
   \emph{strict homomorphism} 
   $\mathcal{F} \colon \bicat{A} \to\bicat{B}$ consists
   of a function 
   $\mathcal{F} \colon \Ob\bicat{A} \to\Ob\bicat{B}$
   and a family of functors
   $\mathcal{F}(A,B) \colon
   \bicat{A}(A,B) \to \bicat{B}(\mathcal{F}A,\mathcal{F}B)$ such that:
   \begin{enumerate}[(i)]
   \item for all $A\in\bicat{A}$,
     $\mathcal{F}(\id_A)=\id_{\mathcal{F}A}$;
   \item for all $f \colon A \to B$, 
     $g \colon B \to C$ in $\bicat{A}$,
     $\mathcal{F}(g \circ f)=\mathcal{F}(g) \circ \mathcal{F}(f)$;
   \item $\mathcal{F}\alpha_{A,B,C,D}
     =\alpha_{\mathcal{F}A,\mathcal{F}B,\mathcal{F}C,\mathcal{F}D}$.
   \end{enumerate}
 \end{defi}

 \begin{exa}\label{exm:homs}
   The following strict homomorphisms will be of interest to us:
   \begin{enumerate}[$\bullet$]
   \item the covariant embedding $\Gr \colon \cat{C} \to\Span[C]$ 
     which acts as the identity
     on objects and takes an arrow $f \colon C \to D$ to its graph
     $\spa{\id}{f} \colon C\rightharpoonup D$. 
     It allows to
     consider the objects and morphisms in the ``universe'' $\cat{C}$
     in the ``larger universe'' $\Span[C]$; 
   \item $\Gamma\diag{F} \colon \cat{J} \to\Span[C]$ where
     $\diag{F} \colon \cat{J} \to \cat{C}$ is a functor.
     It allows to lift every diagram in $\cat{C}$ to a diagram in $\Span[C]$;
   \item given a \textsc{su}-bicategory $\bicat{B}$ and $B\in\Ob\bicat{B}$,
     we shall abuse notation and denote the strict homomorphism 
     from $\cat{J}$ to $\bicat{B}$ 
     which is constant at $B$ by $\Delta[\CAT{J}] B$. 
     It is used to define conical bi-colimits for diagrams.
     Note that in the case of $\bicat{B}=\Span[C]$, 
     ``$\Delta=\Gamma\Delta$''.
   \end{enumerate}
 \end{exa}

   \begin{defi}[Lax transformations]\rm%
     \label{defn:lax}
     Given strict homomorphisms 
     $\mathcal{F},\mathcal{G} \colon \bicat{A} \to \bicat{B}$ between \textsc{su}-bicategories,
     a \emph{(lax) transformation} consists of arrows
     $\kappa_A \colon \mathcal{F}A \to \mathcal{G}A$ for 
     all $A\in\bicat{A}$ and 2-cells 
     $\kappa_f \colon \mathcal{G}f\circ\kappa_A \Rightarrow 
     \kappa_B \circ \mathcal{F}f$ for all $f \colon A \to B$ in $\bicat{A}$ (illustrated below)
   \[
     \begin{tikzpicture}[scale=1.5]
       \mnode[FA]{\mathcal{F}A}{0,0}
       \mnode[FB]{\mathcal{F}B}{1,0}
       \mnode[GA]{\mathcal{G}A}{0,-1}
       \mnode[GB]{\mathcal{G}B}{1,-1}
       \ardrawL{FA}{GA}{left}{\kappa_A}
       \ardrawL{FB}{GB}{right}{\kappa_B}
       \ardrawL{FA}{FB}{above}{\mathcal{F}f}
       \ardrawL{GA}{GB}{below}{\mathcal{G}f}
       \begin{scope}[shorten >=1ex,shorten <=1ex]
         \artwoL{GA}{FB}{above,pos=.3}{\kappa_f~~}   
       \end{scope}
     \end{tikzpicture}
   \]
     such that: 
 \begin{enumerate}[(i)]
 \item $\kappa_{\id_A}=\iota_{\kappa_A}$ for each $A\in\bicat{A}$; 
 \item for any $f \colon A \to B$, $g \colon B \to C$ in $\bicat{A}$, 
   the following 2-cells are equal:
   \begin{displaymath}
     \begin{tikzpicture}[scale=2.3,yscale=1.3,baseline={(0,-1.3)}]
       \mnode[FA]{\mathcal{F}A}{0,0}
       \mnode[FB]{\mathcal{F}B}{1,0}
       \mnode[FC]{\mathcal{F}C}{2,0}
       \ardrawL{FA}{FB}{above}{\mathcal{F}f}
       \ardrawL{FB}{FC}{above}{\mathcal{F}g}
       \mnode[GA]{\mathcal{G}A}{0,-1}
       \mnode[GB]{\mathcal{G}B}{1,-1}
       \mnode[GC]{\mathcal{G}C}{2,-1}
       \ardrawL{GA}{GB}{below}{\mathcal{G}f}
       \ardrawL{GB}{GC}{below}{\mathcal{G}g}
       \ardrawL{FA}{GA}{left}{\kappa_A}
       \ardrawL{FC}{GC}{right}{\kappa_C}
       \draw[->](FA) -- (GB) coordinate[pos=.7](mid);
       \ardrawL{FA}{GB}{above right,pos=.3}{\!\kappa_B \circ \mathcal{F}f} 
       \draw[->](FB) -- (GC)coordinate[pos=.3](mid');
       \ardrawL{FB}{GC}{below left,pos=.7}{\mathcal{G}g \circ \kappa_B\!}
       \begin{scope}[shorten >=.5ex,shorten <=.5ex]
         \artwoL{GA}{mid}{left,pos=.8}{\kappa_f\,}
         \artwoL{[xshift=.5ex]mid}{[xshift=-.5ex]mid'}{above,pos=.2}{\alpha}
         \artwoL{mid'}{FC}{below right}{\!\!\kappa_g}
       \end{scope}
     \end{tikzpicture}
     \quad = \quad 
     \begin{tikzpicture}[scale=2.3,yscale=1.3,baseline={(0,-1.3)}]
       \mnode[FA]{\mathcal{F}A}{0,0}
       \mnode[FB]{\mathcal{F}B}{1,0}
       \mnode[FC]{\mathcal{F}C}{2,0}
       \ardrawL{FA}{FB}{above}{\mathcal{F}f}
       \ardrawL{FB}{FC}{above}{\mathcal{F}g}
       \mnode[GA]{\mathcal{G}A}{0,-1}
       \mnode[GB]{\mathcal{G}B}{1,-1}
       \mnode[GC]{\mathcal{G}C}{2,-1}
       \ardrawL{GA}{GB}{below}{\mathcal{G}f}
       \ardrawL{GB}{GC}{below}{\mathcal{G}g}
       \ardrawL{FA}{GA}{left}{\kappa_A}
       \ardrawL{FC}{GC}{right}{\kappa_C}
       \draw[->](FA) to[bend right=18] (GC) 
       	node[rectangle,fill=white,pos=.2]{\color{black}$\scriptstyle \mathcal{G}( g f )  \kappa_A$} 
       	coordinate[pos=.4](mid);       
       \draw[->](FA) to[bend left=18] (GC) 
	    node[rectangle,fill=white,pos=.8]{\color{black}$\scriptstyle \kappa_C \mathcal{F}(gf)$}        
       coordinate[pos=.6](mid');
       \begin{scope}[shorten >=.5ex,shorten <=.5ex]
         \artwoL{GA}{mid}{left,pos=.8}{\alpha\,\,}
         \artwoL{mid}{mid'}{above,pos=.5}{\kappa_{g \circ f}}
         \artwoL{mid'}{FC}{below right}{\alpha}
       \end{scope}
     \end{tikzpicture}.
   \end{displaymath}
 \end{enumerate}
 A transformation is said to be \emph{strong} when all the $\kappa_f$ are
 invertible 2-cells. 
 Given $B\in\bicat{B}$ and a homomorphism $\diag{M} \colon \cat{J} \to\bicat{B}$,
 a \emph{pseudo-cocone} $\lambda \colon \diag{M} \to \Delta B$
 is a synonym for a
 strong transformation $\lambda \colon \diag{M} \to\Delta B$.
 \end{defi}

 \medskip
 Because bicategories have 2-cells, there are morphisms between
 transformations. 
 They are called \emph{modifications} and are defined as follows.
 \begin{defi}[Modifications~\cite{B'enabou1967,Kelly1973}]\rm%
   \label{defn:modification:VK-span}
     Given natural transformations $\kappa,\lambda$ from $\mathcal{F}$ to $\mathcal{G}$, 
     a modification $\Xi \colon \kappa \to \lambda$ consists of 2-cells
     $\Xi_A \colon \kappa_A \to\lambda_A$ for $A\in\bicat{A}$
     such that, for all $f \colon A \to B$ in~$\bicat{A}$, 
     $\lambda_f\circledcirc(\iota_{\mathcal{G}f}\ast \Xi_A) = 
     (\Xi_B\ast \iota_{\mathcal{F}f})\circledcirc \kappa_f$. 
   \[  
     \begin{tikzpicture}[baseline={(FA.center)},scale=2]
       \mnode[FA]{\mathcal{F}A}{0,0}
       \mnode[FB]{\mathcal{F}B}{1,0}
       \mnode[GA]{\mathcal{G}A}{0,-1}
       \mnode[GB]{\mathcal{G}B}{1,-1}
       \ardrawL[{to[bend right=32]}]{FA}{GA}{left,pos=.7}{\kappa_A}
       \ardrawL[{to[bend left=32]}]{FA}{GA}{right,pos=.3}{\lambda_A}
       \draw[bend right=32,draw=none] (FA) to (GA) coordinate[pos=.7](mid);
       \draw[bend left=32,draw=none] (FA) to (GA) coordinate[pos=.7](mid');
       \artwoL{[xshift=.5ex]mid}{[xshift=-.5ex]mid'}{above}{\Xi_A}
       \begin{scope}[shorten <=1.5ex,shorten >=1.5ex]
         \artwoL[{to[bend right=20]}]{GA.north east}{FB}{below,pos=.65}{~~\,\lambda_f}
       \end{scope}
       \ardrawL[{to[bend left=32]}]{FB}{GB}{right,pos=.3}{\lambda_B}
       \ardrawL{FA}{FB}{above}{\mathcal{F}f}
       \ardrawL{GA}{GB}{below}{\mathcal{G}f}
     \end{tikzpicture}
     \raisebox{-5ex}{$\qquad=\qquad$}
     \begin{tikzpicture}[baseline={(FA.center)},scale=2]
       \mnode[FA]{\mathcal{F}A}{0,0}
       \mnode[FB]{\mathcal{F}B}{1,0}
       \mnode[GA]{\mathcal{G}A}{0,-1}
       \mnode[GB]{\mathcal{G}B}{1,-1}
       \ardrawL[{to[bend right=32]}]{FB}{GB}{left,pos=.7}{\kappa_B}
       \ardrawL[{to[bend left=32]}]{FB}{GB}{right,pos=.3}{\lambda_B}
       \draw[bend right=32,draw=none] (FB) to (GB) coordinate[pos=.3](mid);
       \draw[bend left=32,draw=none] (FB) to (GB) coordinate[pos=.3](mid');
       \artwoL{[xshift=.5ex]mid}{[xshift=-.5ex]mid'}{below}{\Xi_B}
       \begin{scope}[shorten >=1.5ex,shorten <=1.5ex]
         \artwoL[{to[bend left=20]}]{GA}{FB.south west}{above,pos=.3}{\kappa_f~~}
       \end{scope}
       \ardrawL[{to[bend right=32]}]{FA}{GA}{left,pos=.7}{\kappa_A}
       \ardrawL{FA}{FB}{above}{\mathcal{F}f}
       \ardrawL{GA}{GB}{below}{\mathcal{G}f}
     \end{tikzpicture}
   \]
   Composition is component\-wise, 
   the identity modification on~$\kappa$ is
   ${I_\kappa = \{\iota_{\kappa_{A}}\}_{A \in \mathscr{A}}}$.
 \end{defi}
 Given \textsc{su}-bicategories $\bicat{A}$ and
 $\bicat{B}$, let $\Hom_l{[\bicat{A},\bicat{B}]}$ denote the
 \textsc{su}-bicategory of homomorphisms, lax transformations and modifications.
 Let $\Hom{[\bicat{A},\bicat{B}]}$ denote the corresponding 
 \textsc{su}-bicategory with arrows the strong transformations.


\section{Van Kampen cocones}
\label{sec:van-kampen-diagrams}


\noindent The extensivity condition for coproducts and the notion of Van Kampen condition for
pushouts are both examples (for particular colimits) of a more general condition.
Colimit diagrams that satisfy it are called Van Kampen cocones. Here we give 
this definition together with an elementary characterisation.

%
%
Let us consider coproducts as a motivating example.
A coproduct diagram 
\[A \to[{i_1}] A + B \gets[{i_2}] B\]
in a category $\CAT{C}$ is a cocone of the two-object 
diagram $\langle A,B\rangle$.
If $\CAT{C}$ has chosen pullbacks (along coproduct injections) 
then $i_1$ yields a functor $i_1^*\colon\cat{C}\downarrow(A+B)\to\cat{C}\downarrow A$
that takes an arrow $x\colon X\to A+B$ to its pullback ${i_1}^*x \colon {i_1}^*X  \to A$  
along $i_1$---similarly for $i_2$. 
%
Then $x\mapsto \langle{i_1}^*x, {i_2}^*x\rangle$ defines the functor
$\langle{i_1}^*\textvisiblespace,
{i_2}^*\textvisiblespace\rangle \colon 
(\CAT{C}\,{\downarrow}\, A+B) \to 
(\CAT{C}\,{\downarrow}\,A \times \CAT{C}\,{\downarrow}\,B)$ on objects.
The coproduct $A+B$ satisfies the extensivity condition exactly
when this functor is an equivalence of categories (see~\cite{Carboni1993}).

The situation readily generalises as follows:
replace $\langle{i_1} ,{i_2}\rangle$ by any cocone
$\cocone{\kappa} \colon \diag{D} \to \Delta{A}$ 
from a functor $\diag{D} \colon \CAT{J} \to \CAT{C}$ 
to an object $A$
in a category $\CAT{C}$ 
 with (enough) pullbacks.
Any arrow $x \colon X \to A$  
induces a natural transformation
 $\Delta{x} \colon \Delta{X} \to \Delta{A}$ in the obvious way.
 Since also $\cocone{\kappa} \colon \diag{D} \to \Delta{A}$ 
 is a natural transformation, $\Delta{x}$
 can be pulled back along $\kappa$ in 
 the functor category $[\CAT{J},\CAT{C}]$ yielding a
 natural transformation
 $\cocone{\kappa}^*(\Delta x) \colon
 \cocone{\kappa}^*(\Delta X) \to \diag{D}$. 
 \[
 \begin{tikzpicture}[scale=1.5]
   \begin{scope}[xshift=-3em]
     \mnode{A}{0,0}
     \mnode{X}{0,1}
     \ardrawL{X}{A}{left}{x}
   \end{scope}
   \mnode[DA]{\Delta A}{0,0}
   \mnode[DX]{\Delta X}{0,1}
   \ardrawL{DX}{DA}{left}{\Delta x}
   \mnode[D]{\mathcal{D}}{1,0}
   \mnode[X*]{\mkern 20mu \kappa^*(\Delta X)\mkern-20mu}{1,1}
   \ardrawL{X*}{D}{right}{\kappa^*(\Delta x)}
   \begin{scope}[shorten <=-1em]
     \ardrawL{X*}{DX}{ }{}  
   \end{scope}
   \ardrawL{D}{DA}{below}{\kappa}
   \POC[.14]{X*}{-135}
 \end{tikzpicture}
 \]
The described operation extends to a functor
$\cocone{\kappa}^*(\Delta\textvisiblespace)$ 
from $\CAT{C}\,{\downarrow}\, A$
to (a full subcategory of) $[\CAT{J},\CAT{C}] \,{\downarrow}\, \diag{D}$
using the universal property of pullbacks;
it takes morphisms with codomain $A$ to
\emph{cartesian transformations} with codomain $\diag{D}$.
\begin{defi}[Cartesian transformations]\rm%
 \label{def:cartesian-transformation:VK-span}
 Let $\diag{E},\diag{D} \in [\CAT{J}, \CAT{C}]$ be functors
 and let ${\natra{\tau} \colon \diag{E} \to \diag{D}}$ be a
 natural transformation. 
 Then $\natra{\tau}$ is \emph{cartesian} 
 when all the naturality squares are pullback squares, 
 i.e.\ 
 if the pair $\diag{E}_i \gets[\natra{\tau}_i] \diag{D}_i \to[\diag{D}_u]\diag{D}_j$ 
 is a pullback of
 $\diag{E}_i \to[\diag{E}_u]\diag{E}_j \gets[\natra{\tau}_{j}]\diag{D}_j$
 for all $u \colon i \to j$ in $\CAT{J}$.
 \[
 \begin{tikzpicture}[scale=1.5]
     \mnode[Ei]{\diag{E}_i}{0,0}
     \mnode[Ej]{\diag{E}_j}{1,0}
     \ardrawL{Ei}{Ej}{above}{\diag{E}_u}
     \mnode[Di]{\diag{D}_i}{0,-1}
     \mnode[Dj]{\diag{D}_j}{1,-1}
     \ardrawL{Di}{Dj}{below}{\diag{D}_u}
     \ardrawL{Ei}{Di}{left}{\tau_i}
     \ardrawL{Ej}{Dj}{right}{\tau_j}
     \POC[.14]{Ei}{-45}
   \end{tikzpicture}
 \]
   Let $[\CAT{J},\CAT{C}] \,{\downarrow}\, \diag{D}$ 
   be the slice category over $\diag{D}$, 
   which has natural transformations with codomain $\diag{D}$ as objects.
   Let $[\CAT{J},\CAT{C}] \,{\Downarrow}\, \diag{D}$ 
   denote the full subcategory of 
   $[\CAT{J},\CAT{C}] \,{\downarrow}\, \diag{D}$ 
   with the cartesian transformations as objects.
\end{defi}

\begin{defi}[Van Kampen cocones]\rm 
 \label{def:VK-diagram:VK-span}
 Let $\diag{D} \colon \CAT{J} \to \CAT{C}$
 be a functor, and let
 $\cocone{\kappa} \colon \diag{D} \to \Delta A$
 be a cocone such that pullbacks along each 
 $\cocone{\kappa}_i$ exist $(i \in \CAT{J})$. 
 Then $\cocone{\kappa}$ is \emph{Van Kampen \textsc{(vk)}} 
 if the functor 
 $\cocone{\kappa}^*(\Delta\textvisiblespace) 
 \colon \CAT{C}\,{\downarrow}\,A \to
 {[\CAT{J},\CAT{C}]} \,{\Downarrow}\, \diag{D}$
 is an equivalence of categories.
\end{defi}

Extensive and adhesive categories have elementary 
characterisations that are special cases of the following. 
\begin{prop}[Elementary \scalebox{.9}{\small VK} characterisation]
 \label{prop:VK-characterization:VK-span}
 Suppose that $\CAT{C}$ has pullbacks and $\CAT{J}$-colimits,
 $\diag{D} \colon \CAT{J} \to \CAT{C}$
 is a functor and
 $\cocone{\kappa} \colon \diag{D} \to \Delta[J]A$ 
 a cocone
 such that $\CAT{C}$ has pullbacks along 
 $\cocone{\kappa}_i$ ($i \in \CAT{J}$).
 Then 
 $\cocone{\kappa} \colon \diag{D} \to \Delta[J]A$ 
 is Van Kampen 
 iff for every cartesian transformation 
 $\natra{\tau} \colon \diag{E} \to \diag{D}$, 
 arrow $x \colon X \to A$ and cocone
 $\cocone{\beta} \colon \diag{E} \to \Delta[J]X$ 
 such that
 $\natra{\kappa} \circ \cocone{\tau} 
 = \Delta{x} \circ \cocone{\beta}$,
 the following are equivalent:
 \begin{enumerate}[\em(i)]
 \item
   $\cocone{\beta} \colon \diag{E} \to \Delta[J]X$ 
   is a $\CAT{C}$-colimit;
 \item 
   $\diag{D}_i \gets[\natra{\tau}_i]\diag{E}_i \to[\cocone{\beta}_i]X$
   is a pullback of
   $\diag{D}_i \to[\cocone{\kappa}_i]A\gets[x]X$ 
   for all $i \in \CAT{J}$.
 \end{enumerate}
   \[
     \begin{tikzpicture}[scale=1.7]
        \mnode[Ei]{\mathcal{E}_i}{0,0}
        \mnode[Ej]{\mathcal{E}_j}{1,0}
        \mnode[Di]{\mathcal{D}_i}{0,-1}
        \mnode[Dj]{\mathcal{D}_j}{1,-1}
        \ardrawL{Ei}{Ej}{above}{\mathcal{E}_u}
        \ardrawL{Di}{Dj}{above}{\mathcal{D}_u}
        \ardrawL{Ei}{Di}{left,pos=.8}{\tau_i}
        \ardrawL{Ej}{Dj}{left,pos=.8}{\tau_j}
        \mnode[E]{X}{2,0,1}
        \mnode[C]{A}{2,-1,1}
        \ardrawL{Ej}{E}{above right,pos=.3}{\beta_j}
        \ardrawL{Dj}{C}{above right,pos=.3}{\kappa_j}
        \POC[.11]{Ei}{-45}
        \ardrawLd[..controls+(1,0,1)..]{Ei}{E}{below left,pos=.4}{\beta_{i}}
        \ardrawL[..controls+(1,0,1)..]{Di}{C}{below left,pos=.4}{\kappa_{i}}
        \ardrawL{E}{C}{right}{x}
     \end{tikzpicture}
   \]
\end{prop}
\begin{proof}
The proof is a straightforward generalisation of the corresponding characterisation
of Van Kampen squares~\cite[Proposition~2.6]{Lack2005}. 
Assume (i)$\Leftrightarrow$(ii) as well as existence of $\cat{J}$-colimits.
Essential surjectivity of $\kappa^*$ follows easily from (i)$\Rightarrow$(ii).
The fact that $\kappa^*$ is faithful follows from (ii)$\Rightarrow$(i)
and it is full because of existence of $\cat{J}$-colimits.

Conversely, in the presence of $\cat{J}$-colimits, $\kappa^*$ has a left
adjoint defined in the obvious way by taking the appropriate colimits. Then if 
$\kappa^*$ is an equivalence, it follows that the unit and counit of this
adjunction are isomorphisms. Now (i)$\Rightarrow$(ii) is implied
by the fact that the unit is an isomorphism and (ii)$\Rightarrow$(i)
is a consequence of the counit being an isomorphism.
\end{proof}
Cockett and Guo's~\cite{Cockett2007} definition of Van Kampen colimits
is the equivalence of (i) and (ii) in our
Proposition~\ref{prop:VK-characterization:VK-span}.
If the relevant pullbacks and colimits exist then clearly 
the two definitions are equivalent.
\begin{rem}
 With the assumptions of Proposition~\ref{prop:VK-characterization:VK-span},
 any Van Kampen cocone $\cocone{\kappa} \colon \diag{D} \to \Delta A$ 
 is a colimit diagram of $\diag{D}$ in $\cat{C}$ (take $\tau=\id_{\diag{D}}$ and $x=\id_A$).
\end{rem}
\begin{exa}\rm
 \label{exm:VKdiags}
 The following well-known concepts are examples of \vk-cocones:
 \begin{enumerate}[(i)]
 \item a strict initial object is
   a \vk-cocone for the functor from the empty category;
 \item a coproduct diagram in an extensive category~\cite{Carboni1993} is a
   \vk-cocone for a functor from the discrete two object category; 
  \item a regular epimorphism $p \colon E \to B$ in a regular category
    is the Van Kampen cocone \emph{of its kernel pair}%
    \footnote{%
      The kernel pair of a morphism is obtained by pulling it back along itself.
    }%
    \ $r,s \colon P \rightrightarrows E$, 
    i.e.\ in any commuting diagram of the following form
    \begin{displaymath}
      \begin{tikzpicture}[scale=1.5]
        \mnode{Q}{-1,1}
        \mnode{F}{0,1}
        \begin{scope}[shorten > = 2pt]
          \ardrawL{Q.20}{F.160}{above}{u}      
        \end{scope}
        \begin{scope}[shorten < = 2pt]
          \ardrawL{Q.-20}{F.-160}{below}{\,\,v}
        \end{scope}
        \mnode{C}{1,1}
        \ardrawL{F}{C}{above}{q}
        \mnode{P}{-1,0}
        \mnode{E}{0,0}
        \begin{scope}[shorten > = 2pt]
          \ardrawL{P.20}{E.160}{above}{r}
        \end{scope}
        \begin{scope}[shorten < = 2pt]
          \ardrawL{P.-20}{E.-160}{below}{\,\,s}
        \end{scope}
        \mnode{B}{1,0}
        \ardrawL{E}{B}{below}{p}
        \foreach \u/\v/\w in {Q/P/x,F/E/y,C/B/z} 
        {
          \ardrawL{\u}{\v}{left}{\w}
        };
        \POC{Q.-75}{-45}
        \POC{Q}{-45}
      \end{tikzpicture}
    \end{displaymath}
    in which $u,v \colon Q \rightrightarrows F$ is a kernel pair (of some morphism)
    and the squares on the left are pullbacks, 
    the morphism $q$ is the coequaliser of $u,v \colon Q \rightrightarrows F$ 
    if and only if
    the right square is a pullback. This is a direct consequence of a result of 
    Barr and Kock~\cite[Example~6.10]{Barr1971}. Notice that in this case 
    the colimits in question are coequalisers \emph{of kernel pairs}; this
    means that we implicitly restrict the category on the right hand side of the equivalence
    in Definition~\ref{def:VK-diagram:VK-span} to the full subcategory with
    objects the kernel pairs.
 \end{enumerate}
\end{exa}
A \textsc{vk}-cocone from a span  
is a Van Kampen square (see Definition~\ref{def:VK-span-later:diss}).
 In the left hand diagram in Figure~\ref{fig:vanKampen-once-more-span:diss}, 
 the two arrows $B \gets[f] A \to[m] C$ 
 describe a diagram from the three object category
 ${\cdot}\gets{\cdot} \to{\cdot}$, 
 and the cospan $B \to[n]D \gets[g] C$ gives a cocone for this diagram.
 That the back faces are pullback squares
 means that we have a cartesian transformation 
 from $B' \gets[f'] A' \to[m'] C$ to $B \gets[f] A \to[m] C$.
Adhesive categories are thus precisely categories with pullbacks 
in which pushouts along monomorphisms exist and are \textsc{vk}-cocones.


\section{Colimits in the span bicategory}
\label{sec:vkinbicats}


 \noindent In Section~\ref{sec:prelude} we showed that coproducts that satisfy the
 extensivity condition are precisely those that are preserved into the span
 category (Proposition~\ref{prop:characterisationOfExtensivity:vks}) when the latter is
 well-defined, i.e.\ if all pullbacks exist.
 This correspondence failed when we considered pushouts which are Van Kampen 
 squares---a condition for pushouts analogous to the extensivity condition for coproducts.
 The main result of this paper, Theorem~\ref{thm:VKcharacterisation}, asserts that
 the correspondence holds not only for pushouts but for general Van Kampen cocones  
 when considering the canonical embedding into the span \emph{bicategory} providing
 that it exists, i.e.\ the ambient category has all pullbacks. 
 In order to show this we shall
 need to recall the canonical notion of colimit in bicategories in general
 and in bicategories of spans in particular.
 
 To give the flavour of the correspondence
 between \textsc{vk}-cocones in $\cat{C}$ and
 colimits in $\Span[C]$ we outline how  
 Van Kampen squares induce (bi)pushout squares  
 via the embedding~$\Gamma$. 
 An illustration of this 
 is given in Figure~\ref{fig:span-cocone-with-mediating-morphism:diss}.

 \begin{figure}[htb]
   \centering
   \subfloat[{%
     $\SPAN{C}$}-cocone of the $\CAT{C}$-span 
   \mbox{$B \gets[f] A \to[m]C$}
   ]{%
     \label{fig-sub:coconeVKsquare:VK-span}
     \qquad
     \begin{tikzpicture}[baseline=(D.base),scale=2.3]
       \tikzstyle{every node}=[circle,inner sep=0pt]
       \node (D) at (0,0,1) {$\phantom{D}$} ;
       \mnode{B}{-1,0,0} ;
       \node (C) at (1,0,0) {$C\,$} ;
       \node (A) at (0,0,-1) {$A$} ;
       \node (B') at (-1,1,0) {$B'$} ;
       \node (C') at (1,1,0) {$C'$} ;
       \mnode{A'}{-.2,1,-.8};
       \mnode[A'']{A''\!\!}{.2,1,-.8};
       \draw[draw=none] (A') --node{$\scriptstyle \cong$} (A'');
       \node (E) at (0,2,1) {$E$} ;
       \POCS {A'} {-1,0,+1} {0,-1,0}
       \POCS {A''} {+1,0,+1} {0,-1,0}
       \begin{scope}[thick]
         \pardrawL A B {below,near start}{\strt \Gamma f};
       \end{scope}
       \begin{scope}[thick]
         \PardrawL A C {below left,midway}{\Gamma m\,};      
       \end{scope}
       \ardrawL[..controls+(0,-.5,0)..] {A'} {A} {right,pos=.1}{\,a_b} ;
       \ardrawL[..controls+(0,-.5,0)..] {A''} {A} {left,pos=.5}{a_c\,} ;
       \ardrawL{B'} {B} {right,near end}{\,b} ;   
       \ardrawL {C'} {C} {left,near end}{c\,};      
       \ardrawL {A'} {B'} {above,near end}{f'};
       \begin{scope}[shorten <=-1pt]
         \ardrawL {A''} {C'} {above right,pos=.3}{\strt[]m'};
       \end{scope}
       \ardrawL[..controls+(0,.5,0)..] {B'} {E} {left,pos=.1}{{b'}^{\vphantom l}\!} ;
       \ardrawLd[..controls+(0,.5,0)..]{C'} {E} {right,pos=.1}{\,c'} ;   
     \end{tikzpicture}
     \qquad 
   }
   \subfloat[{Mediating $\Span[C]$-morphism from the $\CAT{C}$-cocone
     \mbox{$B \to[n]D\gets[g]C$}}
   ]{%
     \qquad 
     \label{fig-sub:mediatingVKsquare:VK-span}
     \begin{tikzpicture}[baseline=(D.base),scale=2.3]
       \tikzstyle{every node}=[circle,inner sep=0pt]
       \mnode{B}{-1,0,0} ;
       \node (C) at (1,0,0) {$C\,$} ;
       \node (A) at (0,0,-1) {$A$} ;
       \node (B') at (-1,1,0) {$B'$} ;
       \node (C') at (1,1,0) {$C'$} ;
       \mnode{A'}{-.2,1,-.8};
       \mnode[A'']{A''\!\!}{.2,1,-.8};
       \draw[draw=none] (A') --node{$\scriptstyle \cong$} (A'');
       \node (E) at (0,2,1) {$E$} ;
       \begin{scope}[lightgray,thick]
         \POCS {A'} {-1,0,+1} {0,-1,0}
         \POCS {A''} {+1,0,+1} {0,-1,0}
         \begin{scope}[thick]
           \pardrawL A B {below,near start}{\strt \boldsymbol{\Gamma f}}
         \end{scope}
         \begin{scope}[thick]
           \PardrawL A C {below left,midway}{\boldsymbol{\Gamma m}\,}
         \end{scope}
         \ardrawL[..controls+(0,-.5,0)..] {A'} {A} {right,pos=.1}{\,\boldsymbol{a_b} }
         \ardrawL[..controls+(0,-.5,0)..] {A''} {A} {left,pos=.5}{\boldsymbol{a_c}\,} 
         \ardrawL{B'} {B} {right,near end}{\,\boldsymbol b}    
         \ardrawL {C'} {C} {left,near end}{\boldsymbol c\,}      
         \ardrawL {A'} {B'} {above,near end}{\boldsymbol {f'}}
         \begin{scope}[shorten <=-1pt]
           \ardrawL {A''} {C'} {above right,pos=.3}{\strt[]\boldsymbol {m'}};
         \end{scope}
         \ardrawL[..controls+(0,.5,0)..] {B'} {E.west} {left,pos=.1}{\boldsymbol {b'}^{\vphantom l}} ;
         \ardrawLd[..controls+(0,.5,0)..]{C'} {E.east} {right,pos=.1}{\boldsymbol {\,c'}} ;   
       \end{scope}
       \node (D) at (0,0,1) {$D$} ;
       \node (D') at (0,1,1.6) {$D'$} ;
       \mnode[B'']{B''\!}{-1.35,1,.35}
       \mnode[C'']{C''\!}{1.5,1,.5}
       \POCSd {C''} {-1,0,+1} {0,-1,0}
       \POCSd {B''} {+1,0,+1} {0,-1,0}
       \draw[draw=none,sloped] (B'') -- node{$\scriptstyle \cong$} (B');
       \draw[draw=none,sloped] (C'') -- node{$\scriptstyle \cong$} (C');
       \ardrawL[..controls+(0,-.5,0)..]{B''}{B}{left}{d_n}
       \ardrawL[..controls+(0,-.5,0)..]{C''}{C}{right}{d_g}
       \ardrawLd[..controls+(-.5,0,.5)..]{C''}{D'}{below right,pos=.35}{\strt[]g''}
       \ardrawLd[..controls+(.5,0,.5)..]{B''}{D'}{below left,pos=.65}{\strt n''}
       \begin{scope}[shorten >=2pt]
         \ardrawLd[..controls+(0,-.45,0)..] {D'} {D} {right,pos=.35}{d} ;
         \ardrawLd[..controls+(0,+.5,0)..] {D'} {E.south west} {right,pos=.2}{d'} ;
       \end{scope}
       \begin{scope}[thick]
         \PardrawL {B} {D} {below left}{\Gamma n\,};
         \pardrawLd {C} {D} {below,near start}{\Gamma g\strt};
       \end{scope}
     \end{tikzpicture}
     \qquad 
   }
   \caption{Cocones and mediating morphisms consisting of spans}
   \label{fig:span-cocone-with-mediating-morphism:diss}
 \end{figure}

 At the base of
 Figure~\ref{fig:span-cocone-with-mediating-morphism:diss}%
\textsc{\subref{fig-sub:coconeVKsquare:VK-span}}
 is (the image of) a $\CAT{C}$-span $B \gets[f]A \to[m]C$ in $\Span[C]$. 
 Further, if 
 $\spa{b}{b'} \colon B \rightharpoonup E$ 
 and $\spa{c}{c'} \colon C \rightharpoonup E$ 
 are a pseudo-cocone for 
 $B \gets[\Gamma f]A \to[\Gamma m]C$ in $\Span[C]$ 
 then taking pullbacks of $b$ along $f$ and $c$ along $m$ (in~$\CAT{C}$)
 yields isomorphic objects over~$A$, say $a_b$ and $a_c$; 
 as a result we obtain two pullback squares 
 that will be the back faces of a commutative cube. 

 Next, 
 let the bottom of
 Figure~\ref{fig:span-cocone-with-mediating-morphism:diss}%
 \textsc{\subref{fig-sub:mediatingVKsquare:VK-span}}
 be (the image of) a commuting $\CAT{C}$-square, 
 thus yielding another pseudo-cocone of
 $B \leftharpoonup[\Gamma{f}]A\rightharpoonup[\Gamma{m}]C$, 
 namely $B \rightharpoonup[\Gamma{n}]D\leftharpoonup[\Gamma{g}]C$. 
 If there is a mediating morphism $\spa{d}{d'} \colon D \rightharpoonup E$
 to 
 $B \rightharpoonup[\spa{b}{b'}]E \leftharpoonup[\spa{c}{c'}] E$
 of
 Figure~\ref{fig:span-cocone-with-mediating-morphism:diss}%
 \textsc{\subref{fig-sub:coconeVKsquare:VK-span}}
 then pulling back $d$ along $n$ and~$g$ results in morphisms $d_n$ and~$d_g$,
 which are isomorphic to $b$ and~$c$;
 the resulting pullback squares provide the front faces of a cube. 

 Now, 
 if $B \to[n]D\gets[g]C$ is a \textsc{vk}-cocone of $B \gets[f]A \to[m]C$ 
 then such a mediating morphism
 can be constructed by taking $D'$ as the pushout of $B'$ and~$C'$ 
 over either one of $A'$ or $A''$. 
 The morphisms $d \colon D' \to E$ and $d' \colon D' \to D$ 
 arise from the universal property of pushouts, 
 everything commutes and the front faces are pullback squares because of the \textsc{vk}-property. 
 Further this mediating morphism is essentially unique, 
 which means that given any other span
 $\spa{e}{e'} \colon D \rightharpoonup E$ 
 such that both
 $\spa{b}{b'} \cong \spa{e}{e'} \circ \Gamma{n}$ 
 and 
 $\spa{c}{c'} \cong \spa{e}{e'} \circ \Gamma{g}$ hold, 
 the two spans $\spa{e}{e'}$ and $\spa{d}{d'}$ 
 are isomorphic via a unique isomorphism.

 Though this sketch lacks relevant technical details, 
 it gives a good indication
 of the diagrams that are involved
 in the proof of the fact
 that Van Kampen squares in $\CAT{C}$
 induce bipushouts in $\Span[C]$.
 Moreover, also the converse holds, 
 i.e.\ if the image of a pushout is a bipushout 
 in $\Span[C]$ then it is a Van Kampen square.

 \subsection{Span bicolimits}
 \label{sec:spansBicolimits}
 Clearly any diagram in $\Span[C]$ can
 be ``decomposed'' into a diagram in $\cat{C}$:
 each arrow in $\Span[C]$ gives
 two $\CAT{C}$-arrows from a carrier object;
 moreover a 2-cell is a compatible $\CAT{C}$-arrow 
 between the carriers.
 
 We shall start with further observations 
 along these lines.
 Roughly
 we are able to ``drop a dimension'' in the following sense.
 %
 First, it is easy to see that $[\cat{J},\cat{C}]$ inherits a choice of 
 pullbacks from $\cat{C}$.
 In particular, 
 it follows that $\Span[{[\cat{J},\cat{C}]}]$
 is a \textsc{su}-bicategory.
 Now, given $\diag{F},\diag{G}\in[\cat{J},\cat{C}]$
 we note that:
 \begin{enumerate}[$\bullet$]
 \item spans of natural transformations from $\diag{F}$ to $\diag{G}$
   correspond to lax transformations from
   $\Gamma\diag{F}$ to $\Gamma\diag{G}$; and 
 \item morphisms of such spans are the counterpart of modifications.
 \end{enumerate}
 The following lemma makes this precise.
 \begin{lem}%
   \label{lem:FunctorSpans}%
   There is a strict homomorphism
   \[
   \mathbf{\Gamma} \colon 
   \Span[{[\cat{J},\cat{C}]}] \to \Hom_l{[\cat{J},\Span[C]]}
   \]
   that takes $\diag{F} \in [\CAT{J}, \CAT{C}]$ to $\Gamma\diag{F}$ and 
   is full and faithful on both arrows and 2-cells. 
 \end{lem}
 \begin{proof}
   %
   A span of natural transformations
   $\spa{\varphi}{\psi} \colon \mathcal{F} \rightharpoonup \mathcal{G}$
   with carrier $\mathcal{H}$ 
   is mapped to a lax transformation
   from $\Gamma\mathcal{F}$ to $\Gamma\mathcal{G}$
   as follows:
   for each $i \in \CAT{J}$, we put 
   ${\kappa_i := 
   \spa{\varphi_{i}}{\psi_{i}} \colon 
   \mathcal{F}_i \rightharpoonup \mathcal{G}_i}$, 
   and for each morphism $u \colon i \to j$ in $\CAT{J}$, we define a \mbox{2-cell}
   ${\kappa_u \colon
     \Gamma{\mathcal{G}_u} \circ \kappa_i
     \to \kappa_j \circ \Gamma{\mathcal{F}_u}}$
   as sketched below.
   More explicitly, by naturality of $\varphi$ we have that	
   $\mathcal{F}_u \circ \varphi_i 
   = \varphi_j \circ \mathcal{H}_u$ and so 
   the arrow 
   $\kappa_u \colon \diag{H}_i \to \diag{F}_i \times_{\diag{F}_j} \diag{H}_j$
   is the unique one satisfying
   $\varphi_i = \pi_1\circ\kappa_u$ and
   $\mathcal{H}_u = \pi_2\circ\kappa_u$. 
   To check that $\kappa_u$ is a 2-cell it remains verify
   that
   $\psi_j\circ\pi_2\circ\kappa_u
   =\psi_j\circ\diag{H}_u=\diag{G}_u\circ\psi_i$,
   which follows by the naturality of $\psi$.
   \[
       \begin{tikzpicture}[scale=1.5,xscale=1.3]
         \mnode[Fi]{\mathcal{F}_i}{0,0}
         \mnode[Fj]{\mathcal{F}_j}{1,0}
         \ardrawL{Fi}{Fj}{above}{\mathcal{F}_u}
         \begin{scope}[shift={(-.5,-.5)}]
           \mnode[Hi]{\mathcal{H}i}{-.4,-1}     
         \end{scope}
         \mnode[Hj]{\mathcal{H}_j}{1,-1}
         \mnode[Gi]{\mathcal{G}_i}{0,-2}
         \mnode[Gj]{\mathcal{G}_j}{1,-2}
         \ardrawL{Gi}{Gj}{below}{\mathcal{G}_u}
         \ardrawL{Hj}{Gj}{left}{\psi_{j}}
         \ardrawL{Hj}{Fj}{left}{\varphi_{j}}
         \pardrawL[{to[bend left=15]}]%
         {Fj.south east}{Gj.north east}{right}{\kappa_j}
         \ardrawL{Hi}{Gi}{below left}{\psi_i}
         \ardrawL[..controls+(0,1)..]%
         {Hi}{Fi}{above left,pos=.8}{\varphi_i}
         \mnode[X]{\mathcal{F}_i\times_{\mathcal{F}_j}\mathcal{H}_j}{0,-1}     
         \ardrawL{X}{Fi}{left}{\pi_1}
         \ardrawL{X}{Hj}{below}{\pi_2}
         \POC{X.north}{45}
         \begin{scope}[dashed,thick,shorten >=-.5ex]
           \ardrawL{Hi}{X.south west}{above left,pos=.9}{\kappa_u}     
         \end{scope}
         \begin{scope}
           \ardrawL[..controls+(1.5,0)..]{Hi}{Hj}{below,pos=.3}{\mathcal{H}_u}
         \end{scope}
       \end{tikzpicture}
  \]
   Further, a 2-cell between spans
   ${\spa{\varphi}{\psi}, \spa{\varphi'}{\psi'} 
     \colon  \mathcal{F} \rightharpoonup \mathcal{G}}$
   with respective carriers $\mathcal{H},\mathcal{H}'$
   is a natural transformation $\xi \colon \diag{H} \to \diag{H}'$
   satisfying both 
   $\varphi'\circledcirc \xi=\varphi$ and 
   ${\psi'\circledcirc\xi=\psi}$.
   This induces a modification 
   ${\{\xi_{i}\}_{i \in \CAT{J}} \colon
     \mathbf{\Gamma}\spa{\varphi}{\psi} \to
     \mathbf{\Gamma}\spa{\varphi'}{\psi'}
   }$.
   

 It follows from the definition that $\kappa_{\id_i}=\iota_{\diag{H}_i}$. 
 To check the second requirement
 of lax transformations (see Definition~\ref{defn:lax}), 
 consider two arrows $u \colon i \to j$ and $v \colon j \to k$ in $\CAT{J}$.
 Since $\alpha_{\mathcal{F}_u, \kappa_j,\mathcal{G}v}$ 
 and $\alpha_{\kappa_i, \mathcal{G}_u, \mathcal{G}v }$ 
 are identities, 
 one merely has to show
 $\alpha \circledcirc \kappa_{v\circ u} 
 = \kappa_v \ast \iota_{\mathcal{F}_u} \circledcirc\iota_{\mathcal{G}v}\ast \kappa_u$. 
 The latter equation amounts to commutativity of the diagram below---a consequence of functoriality of $\diag{H}$.
 \[
     \begin{tikzpicture}[scale=1.5]
       \mnode[Hi]{\mathcal{H}_i}{0,0}
       \mnode[FiHj]{\diag{F}_i\times_{\diag{F}_j} \diag{H}_j}{-1,-1}
       \mnode[FiHk]{\diag{F}_i\times_{\diag{F}_k}\diag{H}_k}{1,-1}
       \mnode[xx]{\diag{F}_i\times_{\diag{F}_j} 
         (\diag{F}_j\times_{\diag{F}_k} \diag{H}_k)}{0,-2}
       \ardrawL{Hi}{FiHk}{above right}{\kappa_{v \circ u}}
       \ardrawL{FiHk}{xx}{below right,pos=.25}{%
         \alpha_{\mathcal{F}_u,\mathcal{F}v,\kappa_k}}
       \ardrawL{Hi}{FiHj}{above left}{\kappa_u}
       \ardrawL{FiHj}{xx}{below left,pos=.25}{%
         \diag{F}_i\times_{\diag{F}_j} \kappa_v}
     \end{tikzpicture}
 \]
 Faithfulness on arrows is immediate. 
 Conversely, given a lax transformation
 ${\kappa \colon \Gamma\diag{F} \to\Gamma\diag{G}}$, 
 we construct 
 a functor $\diag{H}$ and natural transformations
 $\varphi \colon \diag{H} \to\diag{F}$,
 $\psi \colon \diag{H} \to\diag{G}$ such that
 $\mathbf{\Gamma}\spa{\varphi}{\psi}=\kappa$ as follows:
 let $\diag{H}_i$ be the carrier of the span
 $\kappa_i$ and $\diag{H}_u\Defeq \pi_2\circ\kappa_u$, 
 and further $\varphi_i$ and $\psi_i$ are the left and
 right component of span $\kappa_i$, respectively. 
 Functoriality
 follows directly from the commutativity of the diagram above;
 naturality of
 $\varphi$ and $\psi$ follows from the fact that each $\kappa_u$ is a 2-cell.

 Finally, a 2-cell
 between spans
 ${\spa{\varphi}{\psi}, \spa{\varphi'}{\psi'} \colon
   \mathcal{F} \rightharpoonup \mathcal{G}}$
 with respective carriers $\mathcal{H},\mathcal{H}'$, 
 is a natural transformation
 $\xi \colon \diag{H} \to \diag{H}'$
 satisfying both $\varphi'\circledcirc \xi=\varphi$ 
 and ${\psi'\circledcirc\xi=\psi}$.
 To prove that such a 2-cell induces a modification 
   ${\{\xi_{i}\}_{i \in \CAT{J}} \colon
     \mathbf{\Gamma}\spa{\varphi}{\psi} \to
     \mathbf{\Gamma}\spa{\varphi'}{\psi'}
   }$
 one needs to verify the equality
 $\kappa'_{u} \circledcirc \xi_i = 
 \xi_j \circledcirc \kappa_{u}$; 
 this amounts to the commutativity of the diagram below.
 To see why the diagram commutes, 
 consider the two projections
 $\mathcal{F}_i \gets[\pi_1] 
 \diag{F}_i\times_{\diag{F}_j} \diag{H}'_j 
 \to [\pi_2] \diag{H}'_j$:
 now one has 
 $\pi_1 \circ \kappa'_u \circ\xi_i
 = \varphi_i = 
 \pi_1 \circ (\diag{F}_i\times_{\diag{F}_j} \xi_j) \circ \kappa_u$
 because of $\varphi = \varphi' \circ \xi$   and
 $\pi_2 \circ \kappa'_u \circ\xi_i
 = \mathcal{H}'_u \circ\xi_i 
 = \xi_j \circ \mathcal{H}_{u} 
 = \pi_2 \circ (\diag{F}_i\times_{\diag{F}_j} \xi_j) \circ \kappa_u$ 
 by naturality of $\xi \colon \mathcal{H} \to \mathcal{H}'$. 
 \[
     \begin{tikzpicture}[scale=1.0+(2.0/4)]
       \mnode[Hi]{\mathcal{H}_i}{0,0}
       \mnode[Hi']{\mathcal{H}'i}{-1,-1}    
       \mnode[FiHj]{\diag{F}_i\times_{\diag{F}_j} \diag{H}_j}{1,-1}    
       \mnode[FiHj']{\diag{F}_i\times_{\diag{F}_j} \diag{H}'_j}{0,-2}    
       \ardrawL{Hi}{Hi'}{above left}{\xi_i}
       \ardrawL{Hi'}{FiHj'}{below left,pos=.2}{\kappa'_u}
       \ardrawL{Hi}{FiHj}{above right}{\kappa_u}
       \ardrawL{FiHj}{FiHj'}{below right,pos=.2}{\diag{F}_i\times_{\diag{F}_j} \xi_j}
     \end{tikzpicture}
 \]

 Conversely,
 any modification $\Xi \colon \kappa \to\kappa'$ 
 is a natural transformation
 ${\{\Xi_{i}\}_{i \in \CAT{J}} \colon \diag{H} \to \diag{H}'}$
 where $\mathcal{H}$ and $\mathcal{H}'$ are
 the respective carrier functors---naturality follows directly
 from the commutativity of the diagram above
 (taking $\xi_i = \Xi_i$ and $u = \id_{i}$). 
 \end{proof}
 \begin{cor}
   \label{cor:laxTransSpan:VK-span}
   For any functor $\diag{F}  \in  [\cat{J},\cat{C}]$,
   the strict homomorphism $\mathbf{\Gamma}$
   defines a natural isomorphism between the following two functors
   of type $[\cat{J},\cat{C}] \to\Cat$:
   \[
   \Span[{[\cat{J},\cat{C}]}](\diag{F},\textvisiblespace) \;\cong\;
   \Hom_l{[\cat{J},\Span[C]]}(\Gamma\diag{F},\Gamma \textvisiblespace).
   \]\qed
 \end{cor}

 The above lemma and corollary can be adapted for strong transformations
 instead of lax ones 
 (this will recur when we discuss bicolimits formally). 
 The restriction to strong transformations has a counterpart on the other 
 side of the isomorphism of Corollary~\ref{cor:laxTransSpan:VK-span}:
 we need to restrict to those spans in $\Span[{[\cat{J},\cat{C}]}](\diag{F},\mathcal{G})$
 that have a cartesian transformation from the carrier to $\diag{F}$.

 Recall that a cartesian transformation between functors 
 is a natural transformation with all naturality squares pullbacks 
 (see Definition~\ref{def:cartesian-transformation:VK-span}). 
 It is an easy exercise to show that cartesian transformations 
 include all natural isomorphisms and are closed under pullback.
 Hence---similarly to how one restricts the arrows of a span bicategory to 
 partial maps
 i.e.\ those spans with the left component mono---we let $\CSpan{[\cat{J},\cat{C}]}$ be the (non-full) sub-bicategory of 
 $\Span[{[\cat{J},\cat{C}]}]$ that has as arrows 
 from $\mathcal{F}$ to $\mathcal{G}$ 
 %
 those spans in which the left component 
 is cartesian.
 Adapting the proof of Lemma~\ref{lem:FunctorSpans}, 
 one obtains the following.

 \begin{prop}\label{prop:diag}%
   There is a strict homomorphism
   $\mathbf{\Gamma} \colon 
   \CSpan{[\cat{J},\cat{C}]} \to \Hom{[\cat{J},\Span[C]\makebox[0pt][l]{$]$}}$
   which is full and faithful on both arrows and 2-cells. 
   For any functor $\diag{F} \in  [\cat{J},\cat{C}]$,  
   $\mathbf{\Gamma}$~defines a natural isomorphism 
   between the following functors $[\cat{J},\cat{C}] \to\Cat$:
   \begin{equation*}
     \CSpan{[\cat{J},\cat{C}]}(\diag{F},\textvisiblespace) \;\cong\;
     \Hom{[\cat{J},\Span[C]]}(\Gamma\diag{F},\Gamma \textvisiblespace).
     \eqno{\qEd} 
   \end{equation*}
 \end{prop}
The above lets us pass between diagrams in $\Span[C]$ and
$\cat{C}$: for example the strong transformations of homomorphisms to $\Span[C]$
are those spans of natural transformations of functors to $\cat{C}$ that have
a cartesian first component; the modifications of the former are the morphisms of spans
of the latter.
This observation will be useful when relating the notion of bicolimit in
$\Span[C]$ with the notion of \textsc{vk}-cocone in $\cat{C}$.

 \medskip
 \label{sec:span-bicolimits}
 For our purposes we need to recall only the definition of (conical) bicolimits~\cite{Kelly1989} 
 for functors with domain an (ordinary) small category $\cat{J}$.
 Given a homomorphism
 $\diag{M} \colon \cat{J} \to \bicat{B}$,
 a bicolimit of $\diag{M}$ is an object
 $\bicolim \diag{M}\in\bicat{B}$ together with a pseudo-cocone
 $\kappa \colon \diag{M} \to \Delta(\bicolim \diag{M})$
 such that ``pre-composition'' with $\kappa$
 gives an equivalence of categories
 \begin{equation}\label{eq:bicolimit}
   \bicat{B}(\bicolim \diag{M},X) \simeq
   \Hom{[\cat{J},\bicat{B}]}(\diag{M},
   \Delta X)
 \end{equation}
 that is natural in~$X$ (\ie the right hand side is 
 essentially representable as a functor
 $\lambda X.\Hom{[\cat{J},\bicat{B}]}(\diag{M},
  \Delta X) \colon \bicat{B} \to \Cat$);
 the pair $\tuple{\bicolim \diag{M}, \kappa}$ is referred
 to as \emph{the bicolimit of $\diag{M}$}. 
 We will often speak of
 $\kappa \colon \diag{M}  \to \Delta{\bicolim \diag{M}}$ 
 as a bicolimit  without mentioning
 the pair $\tuple{\bicolim \diag{M}, \kappa}$ explicitly.

 To make the connection with the elementary characterisation 
 of Van Kampen cocones in 
 Proposition~\ref{prop:VK-characterization:VK-span}, 
 we use the fact that equivalences of categories can
 be characterised as full, faithful functors
 that are essentially surjective on objects to derive the
 following equivalent, elementary definition.
 
   \begin{defi}[Bicolimits]\rm%
     \label{def:bicolimit:vks}\label{defn:bicolim}
     Given an \textsc{su}-bicategory $\mathscr{B}$, 
     a category $\CAT{J}$ 
     and a strict homomorphism
     $\diag{M} \colon \cat{J} \to\bicat{B}$, 
     a \emph{bicolimit} for 
     $\diag{M}$ consists of:
     \begin{enumerate}[$\bullet$]
     \item an object $\bicolim \diag{M}\in\bicat{B}$;
     \item a pseudo-cocone
         $\kappa \colon \diag{M} \to \Delta \bicolim \diag{M}$: 
         for each $i\in \cat{J}$
         an arrow
         $\kappa_i \colon \diag{M}_i \to \bicolim \diag{M}$, 
         and for each $u \colon i \to j$ in~$\cat{J}$ an invertible 2-cell
         $\kappa_u \colon \kappa_i \to \kappa_j \circ \diag{M}_u$
         satisfying the axioms 
         required for $\kappa$ to be a strong transformation. 
	\[
	\begin{tikzpicture}[xscale=1.2,scale=1.5,baseline={(Mi.base)}]
           \mnode[Mi]{\mathcal{M}_i}{0,0}
           \mnode[Mj]{\mathcal{M}_j}{1,0}
           \ardrawL{Mi}{Mj}{above}{\mathcal{M}_u\,}
           \mnode[Bic]{\bicolim\mathcal{M}}{0,-1}
           \ardrawL{Mi}{Bic}{left}{\kappa_i}
           \ardrawL[{to[bend left=30]}]{Mj}{Bic}{below right}{\kappa_j}
           \begin{scope}[shorten >=.5ex,shorten <=1.5ex]
             \artwoL{[yshift=1ex]Bic.60}{[yshift=-1ex]Mj.south west}%
             {above,pos=.5}{\kappa_{u}~}
           \end{scope}
     \end{tikzpicture}
     \]
     \end{enumerate}
     The bicolimit satisfies the following universal properties.
     \begin{enumerate}[(i)]
     \item\label{item:1}
         essential surjectivity: \newline
         for any pseudo-cocone
         $\lambda \colon  \diag{M} \to \Delta X$, 
         there exists 
         $h \colon \bicolim \diag{M} \to X$ in~$\bicat{B}$ 
         and an invertible modification
         $\Theta \colon \lambda \to \Delta h\circledcirc\kappa$.
         The pair $\langle h, \Theta\rangle$ 
         is called a \emph{mediating cell} from $\kappa$ to $\lambda$. 
		\[
        \begin{tikzpicture}[scale=1.5,yscale=1.1,baseline={(M.base)}]
             \mnode[M]{\mathcal{M}}{0,0}
             \mnode[X]{\Delta X}{0,-.9}
             \mnode[BM]{\Delta\bicolim{\mathcal{M}}}{2,-.5}
             \ardrawL{M}{X}{left}{\lambda}
             \draw[draw = none] (M) -- (X) coordinate[pos=.7](mid);
             \ardrawL{M}{BM}{above right}{\kappa}
             \begin{scope}[dashed]
               \ardrawL{BM}{X}{below right}{\Delta h}   
               \begin{scope}[shorten <=1.5ex,shorten >=.5ex]
                 \artwoL{mid}{BM.west}{above,pos=.4}{\Theta}    
               \end{scope}
             \end{scope}
           \end{tikzpicture}
        \]
     \item\label{item:2}
         fullness and faithfulness:\newline
         for any $h,h' \colon \bicolim\diag{M} \to X$ in $\bicat{B}$
         and each modification
         $\Xi \colon \Delta h \circledcirc \kappa 
         \to \Delta h' \circledcirc \kappa$,
         there is a unique 2-cell $\xi \colon h \to h'$ 
         satisfying $\Xi = \Delta\xi \ast \iota_{\kappa}$ 
         (and hence $\xi$ is invertible iff $\Xi$ is).
		 \[
         \begin{tikzpicture}[scale=1.5,yscale=.9,baseline={(M.base)}]
           \begin{scope}[yscale=.9]
             \mnode[M]{\mathcal{M}}{0,0}
             \mnode[X]{\Delta X}{0,-2}
             \ardrawL[{to[bend right=18]}]%
             {M}{X}{left,pos=.1}{\Delta h \circledcirc \kappa}
             \draw[draw=none,bend right=18] (M) to (X) coordinate[pos=.6](mid);
             \ardrawL[{to[bend left=18]}]%
             {M}{X}{right,pos=.1}{\Delta h' \circledcirc \kappa}
             \draw[draw=none,bend left=18]%
             (M) to (X) coordinate[pos=.6](mid');
             \artwoL{[xshift=.5ex]mid}{[xshift=-.5ex]mid'}{above}{\Xi}
             \begin{scope}[shift={(2,0)}]
               \mnode[M]{\mathcal{M}}{0,0}
               \mnode[BM]{\Delta\bicolim \mathcal{M}}{0,-.8}
               \draw[draw=none] (BM) --node{\ensuremath{=}}(mid') ;
               \mnode[X]{\Delta X}{0,-2}
               \ardrawL{M}{BM}{right}{\kappa}
               \ardrawL[{to[bend right=40]}]{BM}{X}{left,}{\Delta h}
               \draw[draw=none,bend right=40] (BM) to (X) coordinate[pos=.7](mid);
               \ardrawL[{to[bend left=40]}]{BM}{X}{right,}{\Delta h'}
               \draw[draw=none,bend left=40] (BM) to (X) coordinate[pos=.7](mid');
               \begin{scope}[dashed]
                 \artwoL{[xshift=.9ex]mid}{[xshift=-.7ex]mid'}{above,pos=.4}{\Delta\xi}    
               \end{scope}
             \end{scope}
           \end{scope}
         \end{tikzpicture}
		\]
     \end{enumerate}
   \end{defi}
   Condition~(\ref{item:2}) of this definition
   implies that mediating cells from a bicolimit 
   to a pseudo-cocone are \emph{essentially unique}:
   any two such mediating cells $\langle h,\Theta\rangle$ 
   and $\langle h',\Theta'\rangle$ are isomorphic since
   $\Theta' \circledcirc\Theta^{-1} \colon 
   \Delta h \circledcirc \kappa \to \Delta h' \circledcirc \kappa$ 
   corresponds to a unique invertible 2-cell 
   $\zeta \colon h \to h'$ such that 
   $\Theta' \circledcirc\Theta^{-1} 
   = \Delta\zeta \ast I_{\kappa}$.
 

 To facilitate the exposition of the relationship between
 the bicolimits in $\Span[C]$ and \vk-cocones in $\cat{C}$
 we shall first reformulate the above elementary definition of bicolimits. 
 Given a pseudo-cocone 
 $\kappa \colon \diag{M} \to \Delta C$,
 a morphism $h \colon C \to D$ will be called
 \emph{universal for~$\kappa$} or \emph{\mbox{$\kappa$-}universal} if,
 given any other morphism $h' \colon C \to D$ with a modification
 $\Xi \colon \Delta h\circledcirc \kappa
 \to \Delta h'\circledcirc\kappa$, 
 there exists a \emph{unique} 2-cell $\xi \colon h \to h'$ 
 satisfying $\Xi = \Delta\xi \ast I_{\kappa}$;
 further, a mediating cell $\langle h,\Theta\rangle$ is called universal, 
 if the morphism $h$ is universal.
 The motivation behind this terminology and 
 the slightly redundant 
 statement of the following proposition
 will become apparent in Section~\ref{sec:vanKampenSpans};
 its proof is straightforward.

 \begin{prop}\label{pro:anotherDefinition}
   A pseudo-cocone $\kappa \colon \diag{M} \to \Delta C$ 
   from a diagram $\diag{M}$ to $C$ is
   a bicolimit iff both of the following hold:
   \begin{enumerate}[\em(i)]
   \item for any pseudo cocone $\lambda \colon \diag{M} \to \Delta D$ 
     there is a universal mediating cell 
     ${\langle h \colon C \to D,\, \Theta\rangle}$ from $\kappa$ to $\lambda$;
   \item all arrows $h \colon C \to D$ 
     are universal for $\kappa$.
   \end{enumerate}
 \end{prop}

 We are interested in bicolimits of strict homomorphisms 
 of the form $\Gamma\mathcal{F}$ 
 where $\mathcal{F} \colon \CAT{J} \to \CAT{C}$ 
 is a functor and $\Gamma \colon \CAT{C} \to\Span[C]$
 is the covariant embedding of~$\CAT{C}$.
 The defining equivalence of bicolimits in (\ref{eq:bicolimit}) 
 specialises as follows:
 \begin{align*}
 \Span[C](\bicolim\Gamma\diag{F},X) &\simeq
 \Hom{[\cat{J},{\Span[C]}]}(\Gamma\diag{F},\Delta X).
 \intertext{Using Proposition~\ref{prop:diag}, this is equivalent to:}    
 \Span[C](\bicolim\Gamma\diag{F},X) &\simeq
 \CSpan{[\cat{J},\cat{C}]}(\diag{F},\Delta X).
 \end{align*}
 We shall exploit working in $\CSpan{[\cat{J},\cat{C}]}$ in the following lemma which
 relates the concepts involved in the elementary definition of bicolimits with diagrams
 in~$\cat{C}$. It will serve as the technical backbone of our main theorem.

 \begin{lem}[Mediating cells and universality for spans]%
   \label{lem:correspondence}
     Let $\kappa \colon \diag{F} \to \Delta C$ be a cocone in $\cat{C}$ 
     of a diagram $\diag{F} \in [\CAT{J}, \CAT{C}]$,
     and let $\lambda \colon \Gamma\diag{F} \to \Delta D$ be a pseudo-cocone in $\Span[C]$
     where $\lambda_i=\spa{\varphi_i}{\psi_i}$ for all $i \in \CAT{J}$:
	 \[
     \begin{tikzpicture}[xscale=1,scale=1.1,baseline={(Fi.center)},yscale=.9]
       \mnode[Fi]{\mathcal{F}_i}{0,0}
       \mnode[Fj]{\mathcal{F}_j}{1,0}
       \ardrawL{Fi}{Fj}{above}{\mathcal{F}_u\,}
       \mnode{C}{0,-1}
       \ardrawL{Fi}{C}{left}{\kappa_i}
       \ardrawL[{to[bend left=30]}]{Fj}{C}{below right}{\kappa_j} 
     \end{tikzpicture}
     \qquad
     \begin{tikzpicture}[xscale=1.35,scale=1.1,baseline={(Fi.center)},yscale=.9]
       \mnode[Fi]{\Gamma\mathcal{F}_i}{0,0}
       \mnode[Fj]{\Gamma\mathcal{F}_j}{1,0}
       \pardrawL{Fi}{Fj}{above}{\Gamma\mathcal{F}_u\,}
       \mnode{D}{0,-1}
       \PardrawL{Fi}{D}{left}{\lambda_i}
       \pardrawL[{to[bend left=30]}]{Fj}{D}{below right}{\lambda_j}
       \begin{scope}[shorten >=.5ex,shorten <=1.5ex]
         \artwoL{[yshift=1ex]D.60}{[yshift=-.5ex]Fj.south west}%
         {above,pos=.5}{\lambda_{u}~}
       \end{scope}
     \end{tikzpicture}
    \]
   \begin{enumerate}[\em(i)]
   \item\label{item:mediatingCell}%
       to give a mediating cell
       \[ \tuple{C\xleftarrow{h_1}H\xrightarrow{h_2}D, 
         \Theta \colon  \lambda \to \Delta\spa{h_1}{h_2}\circledcirc \Gamma\kappa} \]
       from $\Gamma\kappa$ to $\lambda$ 
       is to give a cocone
       ${\vartheta \colon \diag{H} \to \Delta H}$ 
       where $\mathcal{H}$ is the carrier functor 
       of the image of $\lambda$ in $\CSpan{[\cat{J},\cat{C}]}(\diag{F},\Delta D)$ 
       (see Proposition~\ref{prop:diag}) such that the resulting
       three-dimensional diagram $(\dagger)$ in $\cat{C}$ (below) commutes 
       and its lateral faces 
       $\pbSqrR{\diag{F}_i}{C}{H}{\diag{H}_i}$ are pullbacks;
	\[
       \begin{tikzpicture}[scale=1.7,baseline={(D.north)}]
         \mnode[Hi]{\mathcal{H}_i}{0,0}
         \mnode[Hj]{\mathcal{H}_j}{1,0}
         \mnode[Fi]{\mathcal{F}_i}{0,-1}
         \mnode[Fj]{\mathcal{F}_j}{1,-1}
         \ardrawL{Hi}{Hj}{above}{\mathcal{H}_u}
         \ardrawL{Fi}{Fj}{above}{\mathcal{F}_u}
         \ardrawL{Hi}{Fi}{left,pos=.8}{\varphi_i}
         \ardrawL{Hj}{Fj}{left,pos=.8}{\varphi_j}
         \mnode{H}{2,0,1}
         \mnode{D}{2,1,1}
         \mnode{C}{2,-1,1}
         \POCSd[.11]{Hj}{0,-1.41,0}{1,0,1}
         \ardrawL{Hj}{H}{above right,pos=.3}{\vartheta_j}
         \ardrawL{Fj}{C}{above right,pos=.3}{\kappa_j}
         \POC[.11]{Hi}{-45}
         \ardrawLd[..controls+(1,0,1)..]{Hi}{H}{below left,pos=.4}{\vartheta_{i}}
         \ardrawL[..controls+(1,0,1)..]{Fi}{C}{below left,pos=.4}{\kappa_{i}}
         \POCSd[.11]{Hi}{0,-1.41,0}{1,0,1}
         \ardrawL{H}{D}{right}{h_2}
         \ardrawL{H}{C}{right}{h_1}
         \ardrawL[..controls+(0,.5,0)..]{Hi}{D}{left,pos=.2}{\psi_i}
         \ardrawL{Hj}{D}{above left,pos=.3}{\psi_j\!}
         \mnode[x]{\qquad(\dagger)}{[xshift=2em]H}
       \end{tikzpicture}
     \]
   \item \label{item:modification}
     to give a modification
     $\Xi \colon
     \Delta\spa{h_1}{h_2}\circledcirc\Gamma\kappa
     \to \Delta\spa{h_1'}{h_2'}\circledcirc\Gamma\kappa$
     for a pair of spans
     \[
     \spa{h_1}{h_2},\spa{h_1'}{h_2'} \colon 
     C\rightharpoonup D
     \]
     is to give a cartesian transformation
     $\Xi \colon \diag{F}\times_{\Delta C} \Delta H 
     \to \diag{F}\times_{\Delta C} \Delta H'$ 
     such that the two equations
     $\pi_1' \circ \Xi=\pi_1$ and
     $(\Delta h_2') \circ \pi_2' \circ \Xi
     =(\Delta h_2) \circ \pi_2$ hold. 
     \begin{align*}
       \begin{tikzpicture}[scale=1.5,baseline={(H.base)},xscale=1.8]
         \mnode[Hi]{\mathcal{F}_i \times_{C} H}{0,0}
         \mnode[Hj]{\mathcal{F}_j\times_{C} H}{1,0}
         \mnode[Fi]{\mathcal{F}_i}{0,-1}
         \mnode[Fj]{\mathcal{F}_j}{1,-1}
         \ardrawL{Hi}{Hj}{above}{\mathcal{F}_u\times_{C} H}
         \ardrawL{Fi}{Fj}{above}{\mathcal{H}_u}
         \ardrawL{Hi}{Fi}{left,pos=.8}{{\pi_1}_i}
         \ardrawL{Hj}{Fj}{left,pos=.8}{{\pi_1}_j}
         \mnode{H}{2,0,1}
         \mnode{D}{2,.7,1}
         \mnode{C}{2,-1,1}
         \POCSd[.11]{Hj}{0,-1.41,0}{1,0,1}
         \ardrawL{Hj}{H}{above right,pos=.3}{{\pi_2}_j}
         \ardrawL{Fj}{C}{above right,pos=.3}{\kappa_j}
         \POC[.11]{Hi}{-45}
         \ardrawLd[..controls+(1,0,1)..]{Hi}{H}{below left,pos=.4}{{\pi_2}_{i}}
         \ardrawL[..controls+(1,0,1)..]{Fi}{C}{below left,pos=.4}{\kappa_{i}}
         \POCSd[.11]{Hi}{0,-1.41,0}{1,0,1}
         \ardrawL{H}{D}{right}{h_2}
         \ardrawL{H}{C}{right}{h_1}
       \end{tikzpicture}
       \qquad
       \hspace{1.5em}
       \begin{tikzpicture}[baseline={(H.base)},scale=.8,xscale=1.1]
         \mnode[Fi]{\mathcal{F}_i}{0,-1}
         \mnode[Hi]{{}\hspace{-1.5em}\mathcal{F}_i\times_{C} H}{-1,0}
         \mnode[H'i]{\mathcal{F}_i\times_{C} H'\hspace{-1.5em}}{1,0}
         \mnode{H}{-1,1}
         \mnode{H'}{1,1}
         \ardrawL{Hi}{H}{left}{{\pi_2}_i}
         \ardrawL{H'i}{H'}{right}{{\pi'_2}_i}
         \mnode{D}{0,2}
         \ardrawL{Hi}{Fi}{below left}{{\pi_1}_i}
         \ardrawL{H'i}{Fi}{below right}{{\pi'_1}_i}
         \ardrawL{H}{D}{above left}{h_2}
         \ardrawL{H'}{D}{above right}{h_2'}
         \ardrawL{Hi}{H'i}{above}{\Xi_i}
       \end{tikzpicture}
       \hspace{1.5em}
       \tag{\ensuremath{\ddagger}}
     \end{align*}
     Here 
     $\mathcal{F} \gets[\pi_1]
     \diag{F}\times_{\Delta C} \Delta H 
     \to[\pi_2] \Delta{H}$ is the pullback of 
     $\mathcal{F} \to[\kappa]\Delta C\gets[\Delta h_1]\Delta{H}$ 
     as sketched in $(\ddagger)$ above
     and similarly for $h' \colon H' \to C$. 
       Further, to give a cell 
       $\xi \colon \spa{h_1}{h_2} \to \spa{h_1'}{h_2'}$ 
       that satisfies 
       $\Delta\xi \ast I_{\Gamma\kappa} = \Xi$ 
       is to give a \mbox{$\CAT{C}$-}arrow
       $\xi \colon H \to H'$ which satisfies the three equations
       $h_1'\circ \xi=h_1$, $h_2'\circ \xi=h_2$ and 
       $\Delta\xi \circ {\pi_2}= {\pi_{2}'} \circ \Xi$;
	\[
       \begin{tikzpicture}[baseline={(D.center)},scale=1.1,xscale=1.4]
         \mnode[Hi]{{}\hspace{-1.5em}\mathcal{F}_i\times_{C} H}{-1,0}
         \mnode[H'i]{\mathcal{F}_i\times_{C} H'\hspace{-1.5em}}{1,0}
         \mnode{H}{-1,1}
         \mnode{H'}{1,1}
         \ardrawL{Hi}{H}{left}{{\pi_2}_i}
         \ardrawL{H'i}{H'}{right}{{\pi'_2}_i}
         \mnode{D}{-1,2}
         \mnode{C}{1,2}
         \ardrawL{H}{D}{left,pos=.7}{h_2}
         \ardrawL[{to[bend right=15]}]{H'}{D}{above right,pos=.95}{h_2'}
         \ardrawL[{to[bend left=15]}]{H}{C}{above left,pos=.95}{h_1}
         \ardrawd[{to[bend left=15]}]{H}{C}
         \ardrawL{H'}{C}{right,pos=.7}{h'_1}
         \ardrawL{Hi}{H'i}{above}{\Xi_i}
         \ardrawL{H}{H'}{above}{\xi}
       \end{tikzpicture}
	\]
   \item\label{item:colimitThenUniversal}
     given a span $\spa{h_1}{h_2} \colon C\rightharpoonup D$,
     if the pullback of $\kappa$ along $h_1$ is a colimit, 
     i.e.\ if $\pi_2 \colon \diag{F}\times_{\Delta C} \Delta H \to \Delta H$
     is a colimit, 
     then $\spa{h_1}{h_2}$ is universal for $\Gamma\kappa$;
   \item \label{item:universalThencolimit}
     conversely, 
     if $\spa{h_1}{h_2}$ is universal for $\Gamma\kappa$, then
     $\pi_2 \colon \diag{F}\times_{\Delta C} \Delta{H} \to \Delta H$ is a 
     colimit---provided that some colimit of
     $\mathcal{F} \times_{\Delta C} \Delta H$ exists in $\CAT{C}$.
   \end{enumerate}
 \end{lem}
 \begin{proof}
   Both (\ref{item:mediatingCell}) and (\ref{item:modification}) 
   are immediate consequences of Proposition~\ref{prop:diag}.

   As for (\ref{item:colimitThenUniversal}), 
   we need to show that every modification
   $\Xi \colon
   \Delta\spa{h_1}{h_2}\circledcirc\Gamma\kappa 
   \to \Delta\spa{h_1'}{h_2'}\circledcirc\Gamma\kappa$ 
   is equal to $\Delta\xi \ast I_{\Gamma\kappa}$ for a unique
   $\xi \colon \spa{h_1}{h_2} \to\spa{h_1'}{h_2'}$. 
   By (\ref{item:modification}), 
   $\Xi$~is a natural transformation
   $\Xi \colon \diag{F}\times_{\Delta C} \Delta{H} 
   \to \diag{F}\times_{\Delta C} \Delta{H'}$.
   Then, by naturality of~$\Xi$, 
   we have that 
   $\pi_{2j}' \circ \Xi_j \circ (\diag{F}_u\times_{\Delta C} \Delta{H}) =
   \pi_{2i}' \circ\Xi_i$ 
   holds for all $u \colon i \to j$ in $\cat{J}$,
   and since $\pi_2$ is a colimit we have a unique
   $\xi \colon H \to H'$ satisfying
   $\xi\circ \pi_{2i}=\pi_{2i}' \circ \Xi_i$ for all $i\in\cat{J}$.
   The equations $h_i = h_i' \circ\xi$ follow from the universal property of $\pi_2$
   (and the properties of $\Xi$). 
   To show uniqueness of $\xi$, 
   let $\zeta \colon \spa{h_1}{h_2} \to \spa{h_1'}{h_2'}$ be a 2-cell 
   such that $\Xi = \Delta\zeta \ast I_{\Gamma_{\kappa}}$;
   then using the second statement of
   Lemma~\ref{lem:correspondence}(\ref{item:modification}), 
   $\Delta\zeta\circ \pi_2 = \pi_2' \circ \Xi$; 
   hence $\zeta=\xi$ follows since $\pi_2$ is a colimit.
   In summary, $\spa{h_1}{h_2}$ is universal for $\Gamma\kappa$.

   To show (\ref{item:universalThencolimit}),
   let $\langle H',\vartheta\rangle$ be a colimit of
   $\diag{F}\times_{\Delta C} \Delta H$. 
   Now, it suffices to show that there is a $\CAT{C}$-morphism
   $\xi \colon H \to H'$ such that
   $\vartheta = \Delta\xi \circ \pi_2$.\footnote{%
     The reason is that once such a $\xi$ is provided, 
     there is a unique $k \colon H' \to H$ 
     satisfying $\Delta k \circ \vartheta=\pi_{2}$, 
     and thus $\xi \circ k = \id_{H'}$ by the universal property of colimits;
     moreover $k \circ \xi=\id_H$ must hold
     since $\spa{h_1}{h_2}$ is universal for $\Gamma\kappa$.}
   By the universal property of $\vartheta$,
   we obtain unique $\CAT{C}$-arrows
   $h_1' \colon H' \to C$ and $h_2' \colon H' \to D$
   such that
   $\Delta{h_1'} \circ \vartheta= \kappa \circ \pi_{1}$ 
   and $\Delta{h_2'} \circ \vartheta= \Delta h_2 \circ \pi_{2}$. 
   It also follows 
   that the two equations $h_1\circ k=h_1'$ and $h_2\circ k = h_2'$ hold.
   Pulling back $\kappa$ along $h_1'$
   yields a span
   $\mathcal{F} \gets[\pi_1']
   \diag{F}\times_{\Delta C} \Delta H'
   \to[\pi_2'] \Delta H'$; 
   we then obtain a natural transformation
   $\Xi \colon 
   \diag{F}\times_{\Delta C} \Delta{H} 
   \to \diag{F}\times_{\Delta C} \Delta{H'}$ 
   which satisfies $\pi_1 = \pi_1' \circ \Xi$ 
   and $\vartheta = \pi_2' \circ \Xi$, 
   and hence also
   $\Delta h_2 \circ \pi_{2}= \Delta{h_2'} \circ \vartheta
   = \Delta{h_2'} \circ \pi_2' \circ \Xi$. 
   By (\ref{item:modification}), this defines a modification
   $\Xi \colon 
   \Delta\spa{h_1}{h_2} \circledcirc \Gamma\kappa 
   \to \Delta\spa{h_1'}{h_2'}\circledcirc\kappa$.
   Using universality, 
   we get a unique $\xi \colon  H \to H'$ such that 
   $h_1' \circ \xi=h_1$, 
   $h_2' \circ \xi=h_2$ and 
   $\Delta \xi \circ \pi_{2} = \pi_{2}' \circ \Xi = \vartheta$.
 \end{proof}


\section{Van Kampen cocones as span bicolimits }
\label{sec:vanKampenSpans}

\noindent Here we prove the main result of this paper,
Theorem~\ref{thm:VKcharacterisation}.
Roughly speaking, the conclusion is that (under natural assumptions---existence of pullbacks and enough colimits in $\CAT{C}$) to be 
\textsc{vk} in $\cat{C}$ is to be a bicolimit in $\Span[\cat{C}]$.
The consequence is that  
``being \textsc{vk}'' is a universal property---in $\Span[C]$ rather than in~$\cat{C}$.

The proof relies on a correspondence 
between the elementary characterisation of
Van Kampen cocones in $\CAT{C}$ of 
Proposition~\ref{prop:VK-characterization:VK-span} 
and  the universal properties of pseudo-cocones in $\Span[C]$ 
of Proposition~\ref{pro:anotherDefinition}.
More precisely, given a colimit  $\kappa \colon \diag{M} \to \Delta C$ in $\CAT{C}$, we shall 
show that:
\begin{enumerate}[$\bullet$]
\item 
 $\Gamma\kappa$-universality \emph{of all} spans
 $\spa{h_1}{h_2} \colon C \rightharpoonup D$
 corresponds to the implication  ${\text{(ii)}\Rightarrow \text{(i)}}$ of Proposition~\ref{prop:VK-characterization:VK-span}, which  is also known as pullback-stability or \emph{universality} of the colimit $\kappa$;
\item 
 \emph{existence of some} 
 universal mediating cell from $\Gamma\kappa$
 to any $\lambda \colon \Gamma\kappa \to \Delta D$ is the counterpart of 
 the implication  ${\text{(i)}\Rightarrow \text{(ii)}}$ of Proposition~\ref{prop:VK-characterization:VK-span}, 
 which---for want of a better name---we here refer to as ``\emph{converse universality}'' of $\kappa$;
\item 
 thus,
 $\Gamma\kappa$ is a bicolimit in $\Span[C]$ 
 if and only if  the colimit $\kappa$ is Van Kampen. 
\end{enumerate}
The first two points are made precise by the statements of the
following two lemmas. The third point is the statement
of the main theorem.

\begin{lem}[Converse universality]%
 \label{lem:converseUniversality}
 Let $\diag{F} \in [\CAT{J}, \CAT{C}]$ 
 where $\CAT{C}$~has pullbacks 
 and for all 
 $(\tau \colon \diag{E} { \to} \diag{F}) \in
 [\CAT{J},\CAT{C}] \,{\Downarrow}\, \diag{F}$
 a colimit of $\diag{E}$ exists. 
 Then 
 $\kappa \colon \diag{F} \to \Delta[J]C$ 
 satisfies ``converse universality'' 
 iff given any pseudo-cocone
 $\lambda \colon \Gamma\diag{F} \to \Delta D$,
 there exists a universal mediating cell
 $\langle\spa{h_1}{h_2},\Theta\rangle$ from
 $\Gamma\kappa$ to $\lambda$ in $\Span[C]$.
\end{lem}
\begin{proof}
 ($\Rightarrow$)
 Suppose that $\lambda \colon \Gamma \diag{F} \to \Delta D$ 
 is a pseudo-cocone in $\Span[C]$.~\smallskip\\
   For $u \colon i \to j$ in $\CAT{J}$, 
   we obtain a commutative diagram, as illustrated
   (see Proposition~\ref{prop:diag}).
   Let $\vartheta \colon \diag{H} \to \Delta H$ 
   be the colimit of~$\diag{H}$;
   thus we obtain $h_1 \colon H \to C$ 
   and $h_2 \colon H \to D$ making diagram~$(\dagger)$ commute.
   By converse universality, the side faces
   $\pbSqrR{\diag{F}_i}{C}{H}{\diag{H}_i}$ are pullback squares;
   using Lemma~\ref{lem:correspondence}(\ref{item:modification})
   we get an invertible modification
   $\Theta \colon \lambda \to
   \Delta\spa{h_1}{h_2} \circledcirc \Gamma\kappa$. 
   That $\spa{h_1}{h_2}$ is universal 
   follows from 
   Lemma~\ref{lem:correspondence}(\ref{item:colimitThenUniversal})
   since $\vartheta$ is a colimit.
\[
   \begin{tikzpicture}[scale=1.5,baseline={(D.base)}]
     \mnode[Hi]{\mathcal{H}_i}{0,0}
     \mnode[Hj]{\mathcal{H}_j}{1,0}
     \mnode[Fi]{\mathcal{F}_i}{0,-1}
     \mnode[Fj]{\mathcal{F}_j}{1,-1}
     \ardrawL{Hi}{Hj}{above}{\mathcal{H}_u}
     \ardrawL{Fi}{Fj}{above}{\mathcal{F}_u}
     \ardrawL{Hi}{Fi}{left,pos=.8}{\varphi_i}
     \ardrawL{Hj}{Fj}{left,pos=.8}{\varphi_j}
     \mnode{D}{2,1,1}
     \mnode{C}{2,-1,1}
     \ardrawL{Fj}{C}{above right,pos=.3}{\kappa_j}
     \POC[.11]{Hi}{-45}
     \ardrawL[..controls+(1,0,1)..]{Fi}{C}{below left,pos=.4}{\kappa_{i}}
     \ardrawL[..controls+(0,.5,0)..]{Hi}{D}{left,pos=.2}{\psi_i}
     \ardrawL{Hj}{D}{above left,pos=.3}{\psi_j\!}
   \end{tikzpicture}
 \] 
 \noindent
   ($\Leftarrow$)
   If in diagram $(\dagger)$ $D=\colim \diag{H}$ 
   and $\langle D,\psi\rangle$ is the corresponding colimit,
   we first use the assumption to obtain a universal 
   mediating cell $\langle\spa{h_1}{h_2}, \Theta\rangle$ 
   from $\Gamma\kappa$ to $\lambda^{(\varphi,\psi)}$ 
   where $\lambda^{(\varphi,\psi)}$ is the pseudo-cocone 
   corresponding
   to the cartesian transformations
   $\varphi \colon \diag{H} \to \diag{F}$ 
   and $\psi \colon \diag{H} \to \Delta D$ such that
   $\lambda_{i}^{(\varphi,\psi)} = \spa{\varphi_i}{\psi_i}$ 
   as in Lemma~\ref{lem:correspondence}(\ref{item:mediatingCell});
   the latter also provides $\vartheta \colon \diag{H} \to \Delta H$ 
   such that $h_2 \circ \vartheta_i =\psi_i$ and all
   $\pbSqrR{\diag{F}_i}{C}{H}{\diag{H}_i}$ are pullback squares.
   \[
   \begin{tikzpicture}[scale=1.8,baseline={(D.base)}]
     \mnode[Hi]{\mathcal{H}_i}{0,0}
     \mnode[Hj]{\mathcal{H}_j}{1,0}
     \mnode[Fi]{\mathcal{F}_i}{0,-1}
     \mnode[Fj]{\mathcal{F}_j}{1,-1}
     \ardrawL{Hi}{Hj}{above}{\mathcal{H}_u}
     \ardrawL{Fi}{Fj}{above}{\mathcal{F}_u}
     \ardrawL{Hi}{Fi}{left,pos=.8}{\varphi_i}
     \ardrawL{Hj}{Fj}{left,pos=.8}{\varphi_j}
     \mnode{H}{2,0,1}
     \mnode{D}{2,1,1}
     \mnode{C}{2,-1,1}
     \POCSd[.11]{Hj}{0,-1.41,0}{1,0,1}
     \ardrawL{Hj}{H}{above right,pos=.3}{\vartheta_j}
     \ardrawL{Fj}{C}{above right,pos=.3}{\kappa_j}
     \POC[.11]{Hi}{-45}
     \ardrawLd[..controls+(1,0,1)..]{Hi}{H}{below left,pos=.4}{\vartheta_{i}}
     \ardrawL[..controls+(1,0,1)..]{Fi}{C}{below left,pos=.4}{\kappa_{i}}
     \POCSd[.11]{Hi}{0,-1.41,0}{1,0,1}
     \ardrawL{H}{D}{right}{h_2}
     \ardrawL{H}{C}{right}{h_1}
     \ardrawL[..controls+(0,.5,0)..]{Hi}{D}{left,pos=.2}{\psi_i}
     \ardrawL{Hj}{D}{above left,pos=.3}{\psi_j\!}
     \mnode[x]{\qquad (\dagger)}{[xshift=1em]H}
   \end{tikzpicture}
\]
 It suffices to show that $h_2 = \id_{H}$.
 However, by the universal property of the colimit
 $\langle{D},{\psi}\rangle$,
 there is an arrow $k \colon D \to h$ 
 such that $k \circ \psi_i=\vartheta_i$.
 The equation $h_2 \circ k=\id_D$ holds because
 $\langle{D},{\psi}\rangle$
 is a colimit in $\CAT{C}$, 
 and $k \circ h_2=\id_{H}$ follows since $\spa{h_1}{h_2}$ 
 is universal for $\Gamma\kappa$.
\end{proof}

\begin{lem}[Universality]\label{lem:universality}
 Consider $\diag{F} \in [\CAT{J}, \CAT{C}]$ 
 where $\CAT{C}$ has pullbacks 
 such that for all
 $(\tau \colon \diag{E} { \to} \diag{F}) 
 \in [\CAT{J},\CAT{C}] \,{\Downarrow}\, \diag{F}$,
 a colimit of $\diag{E}$ exists. 
 Then $\kappa \colon \diag{F} \to \Delta C$ satisfies universality
 iff every morphism
 $\spa{h_1}{h_2} \colon C\rightharpoonup D$ 
 in $\Span[C]$ is universal for $\Gamma\kappa$.
\end{lem}
\begin{proof}
 ($\Rightarrow$)
 Any morphism $\spa{h_1}{h_2}$ leads to a diagram $(\ddagger)$ where 
 all the side-faces are pullbacks. By universality of $\kappa$, 
 the cocone $\pi_2$ of the top face is a colimit;
 thus $\spa{h_1}{h_2}$ is universal for $\Gamma\kappa$ 
 by Lemma~\ref{lem:correspondence}(\ref{item:colimitThenUniversal}).

 \noindent
 ($\Leftarrow$)
 Suppose that in diagram $(\ddagger)$ the side faces are all pullbacks. 
 By assumption $\spa{h_1}{h_2}$ is universal for
 $\Gamma\kappa$, thus
 $\langle H,\pi_2 \colon
 \diag{F}\times_{\Delta C} \Delta H 
 \to \Delta H\rangle$ 
 is a colimit by
 Lemma~\ref{lem:correspondence}(\ref{item:universalThencolimit}). 
\end{proof}

Finally, these two lemmas together with
Proposition~\ref{pro:anotherDefinition}
imply our main result.
\begin{thm}\label{thm:VKcharacterisation}
 Let $\diag{F} \in [\CAT{J}, \CAT{C}]$
 where $\CAT{C}$ has pullbacks 
 and for all cartesian transformations
 $\tau \colon \diag{E} { \to} \diag{F}$,
 a colimit of $\diag{E}$ exists. 
 Then a cocone $\kappa \colon \diag{F} \to \Delta C$ is Van Kampen 
     iff
     $\Gamma\kappa \colon \Gamma\mathcal{F} \to \Delta C$ is a bicolimit
     in $\Span[\cat{C}]$.
 \qed
\end{thm}


\section{Conclusion, related work and future work}
\label{sec:conclusion}

\noindent We gave a general definition of Van Kampen cocone
that captures several previously studied notions
in computer science, topology, and related areas, showing 
that they are instances of the same concept.
Moreover, we have provided two alternative 
characterisations: the first one is elementary, 
and involves only basic category theoretic notions;
the second one exhibits it as a \emph{universal property}:
Van Kampen cocones are just those colimits
that are preserved by the canonical covariant embedding into the span bicategory. 


There is some interesting related recent work.
Milius~\cite{Milius2003} showed that coproducts are preserved (as a
lax-adjoint-cooplimit) in the 2-category of relations over an
extensive category~$\cat{C}$.  Cockett and Guo~\cite{Cockett2007} have
investigated the general conditions under which partial map categories
are join-restriction categories: roughly, certain colimits in the
underlying category are required to be \vk-cocones.

Finally, the definition of Van Kampen cocone allows for several
natural variations.  For example, one may replace the slice category
over the object at the ``tip'' of cocones by a (full) subcategory of
it; this is exactly the step from global descent to
$\mathbb{E}$-descent~\cite{Janelidze1994} and is closely related to
the proposals in~\cite{Ehrig2004a,Ehrig2006a} for a weakening of the
notion of adhesivity.  Alternatively, one may start with cocones or
diagrams of a particular form.  In this way quasi-adhesive
categories~\cite{Lack2005} arise as in the latter only pushouts along
regular monos are required to be \textsc{vk}; another example is the
work of Cockett and Guo~\cite{Cockett2007}, where Van Kampen cocones
exist for a class of diagrams that naturally arises in their study of
join restriction categories.
Thus, possibly combining the latter two ideas,
several new forms of Van Kampen cocones and diagrams arise
as the subject for future research.





%



\begin{thebibliography}{10}

\bibitem{Abramsky1993}
S.~Abramsky.
\newblock Interaction categories.
\newblock In {\em Theory and Formal Methods `93}, Workshops in Computing, pages
  57--69. Springer, 1993.

\bibitem{Aczel1989}
P.~Aczel and N.~Mendler.
\newblock A final coalgebra theorem.
\newblock In {\em Category Theory and Computer Science {(CTCS `89)}}, volume
  389 of {\em LNCS}, pages 357--365. Springer, 1989.

\bibitem{Barr1971}
M.~Barr.
\newblock Exact categories.
\newblock In {\em Exact Categories and Categories of Sheaves}, volume 236 of
  {\em LNM}, pages 1--120. Springer, 1971.

\bibitem{B'enabou1967}
J.~B\'{e}nabou.
\newblock Introduction to bicategories, part 1.
\newblock In {\em Midwest Category Seminar}, volume~47 of {\em LNM}, pages
  1--77. Springer, 1967.

\bibitem{Borceux1994}
F.~Borceux.
\newblock {\em Handbook of categorical algebra}, volume~1.
\newblock Cambridge University Press, 1994.

\bibitem{Brown1997}
R.~Brown and G.~Janelidze.
\newblock {V}an {K}ampen theorems for categories of covering morphisms in
  lextensive categories.
\newblock {\em J.\ Pure Appl.\ Algebra}, 119:255--263, 1997.

\bibitem{Carboni1993}
A.~Carboni, S.~Lack, and R.~F.~C. Walters.
\newblock Introduction to extensive and distributive categories.
\newblock {\em J.\ Pure Appl.\ Algebra}, 84:145--158, February 1993.

\bibitem{Cockett1997}
J.~R.~B. Cockett and D.~A. Spooner.
\newblock Constructing process categories.
\newblock {\em Theor.\ Comput.\ Sci.}, 177(1):73--109, 1997.

\bibitem{Cockett2007}
R.~Cockett and X.~Guo.
\newblock Join restriction categories and the importance of being adhesive.
\newblock Unpublished manuscript, 2007.

\bibitem{Corradini1997}
A.~Corradini, U.~Montanari, F.~Rossi, H.~Ehrig, R.~Heckel, and M.~L\"{o}we.
\newblock Algebraic approaches to graph transformation -- part {\sc i}: Basic
  concepts and double pushout approach.
\newblock In {\em Handbook of Graph Grammars}, pages 163--246. World
  Scientific, 1997.

\bibitem{Ehrig2004a}
H.~Ehrig, A.~Habel, J.~Padberg, and U.~Prange.
\newblock Adhesive high-level replacement categories and systems.
\newblock In {\em Graph Transformation ({ICGT `04})}, volume 3256 of {\em
  LNCS}, pages 144--160. Springer, 2004.

\bibitem{Ehrig2004}
H.~Ehrig and B.~K\"{o}nig.
\newblock Deriving bisimulation congruences in the {DPO} approach to graph
  rewriting.
\newblock In {\em Foundations of Software Science and Computation Structures
  ({FoSSaCS `04})}, LNCS. Springer, 2004.

\bibitem{Ehrig2006a}
H.~Ehrig and U.~Prange.
\newblock Weak adhesive high-level replacement categories and systems: a
  unifying framework for graph and {P}etri net transformations.
\newblock In {\em Essays dedicated to Joseph A. Goguen}, volume 2987 of {\em
  LNCS}, pages 235--251. Springer, 2006.

\bibitem{Freyd1990}
P.~J. Freyd and A.~Scedrov.
\newblock {\em Categories, allegories}.
\newblock North-Holland, 1990.

\bibitem{Gadducci1999}
F.~Gadducci, R.~Heckel, and M.~Llabr\'{e}s.
\newblock A bi-categorical axiomatisation of concurrent graph rewriting.
\newblock In {\em Category Theory and Computer Science ({CTCS `99})}, volume~29
  of {\em ENTCS}. Elsevier, 1999.

\bibitem{Heindel2010a}
T.~Heindel.
\newblock Hereditary pushouts reconsidered.
\newblock In {\em Graph Transformation (ICGT~`10)}, volume 6372 of {\em LNCS},
  pages 250--265. Springer, 2010.
\newblock To appear.

\bibitem{Heindel2009}
T.~Heindel and P.~Soboci\'{n}ski.
\newblock Van {K}ampen colimits as bicolimits in {S}pan.
\newblock In {\em Algebra and Coalgebra in Computer Science ({CALCO `09})},
  number 5728 in LNCS, pages 335--349. Springer, 2009.

\bibitem{Jacobs1995}
B.~Jacobs.
\newblock Parameters and parametrization in specification, using distributive
  categories.
\newblock {\em Fundam.\ Inform.}, 24(3):209--250, 1995.

\bibitem{Janelidze1994}
G.~Janelidze and W.~Tholen.
\newblock Facets of descent, {I}.
\newblock {\em Appl.\ Categor.\ Struct.}, 2(3):245--281, 1994.

\bibitem{Johnstone1977}
P.~T. Johnstone.
\newblock {\em Topos theory}, volume~10 of {\em L.\ M.\ S.\ Monographs}.
\newblock Academic Press, 1977.

\bibitem{Johnstone2002}
P.~T. Johnstone.
\newblock {\em Sketches of an elephant: a topos theory compendium}, volume~1.
\newblock Clarendon Press, 2002.

\bibitem{Johnstone2007}
P.~T. Johnstone, S.~Lack, and P.~Soboci\'{n}ski.
\newblock Quasitoposes, quasiadhesive categories and {A}rtin glueing.
\newblock In {\em Algebra and Coalgebra in Computer Science ({CALCO `07})},
  volume 4626 of {\em LNCS}, pages 312--326. Springer, 2007.

\bibitem{Joyal1994}
A.~Joyal, M.~Nielsen, and G.~Winsk.
\newblock Bisimulation from open maps.
\newblock {\em Inf.\ Comput.}, 127:164--185, 1994.

\bibitem{Katis1997a}
P.~Katis, N.~Sabadini, and R.~F.~C. Walters.
\newblock {Span(Graph):} an algebra of transition systems.
\newblock In {\em Algebraic Methodology and Software Technology (AMAST '97)},
  volume 1349 of {\em LNCS}, pages 322--336. Springer, 1997.

\bibitem{Kelly1989}
G.~M. Kelly.
\newblock Elementary observations on 2-categorical limits.
\newblock {\em Bull.\ Austral.\ Math.\ Soc.}, 39:301--317, 1989.

\bibitem{Kelly1973}
G.~M. Kelly and R.~H. Street.
\newblock Review of the elements of 2-categories.
\newblock In {\em Category Theory}, number 420 in LNM, pages 75--103. Springer,
  1973.

\bibitem{Lack2010}
S.~Lack.
\newblock A 2-categories companion.
\newblock In {\em Towards Higher Categories}. Institute for Mathematics and its
  Applications, 2010.
\newblock To appear.

\bibitem{Lack2004}
S.~Lack and P.~Soboci\'{n}ski.
\newblock Adhesive categories.
\newblock In {\em Foundations of Software Science and Computation Structures
  ({FoSSaCS `04})}, volume 2987 of {\em LNCS}, pages 273--288. Springer, 2004.

\bibitem{Lack2005}
S.~Lack and P.~Soboci\'{n}ski.
\newblock Adhesive and quasiadhesive categories.
\newblock {\em RAIRO-Theor.\ Inf.\ Appl.}, 39(3):511--546, 2005.

\bibitem{Lack2006}
S.~Lack and P.~Soboci\'{n}ski.
\newblock Toposes are adhesive.
\newblock In {\em Graph Transformation ({ICGT `06})}, volume 4178 of {\em
  LNCS}, pages 184--198. Springer, 2006.

\bibitem{Lawvere1991}
F.~W. Lawvere.
\newblock Some thoughts on the future of category theory.
\newblock In {\em Category Theory ({CT `91})}, volume 1488 of {\em LNM}.
  Springer, 1991.

\bibitem{Lindner1976}
H.~Lindner.
\newblock A remark on {M}ackey-functors.
\newblock {\em Manuscripta Math.}, 18(3):273--278, 1976.

\bibitem{Milius2003}
S.~Milius.
\newblock On colimits in categories of relations.
\newblock {\em Appl.\ Categor.\ Struct.}, 11(3):287--312, 2003.

\bibitem{Sassone2005a}
V.~Sassone and P.~Soboci\'{n}ski.
\newblock Reactive systems over cospans.
\newblock In {\em Logic in Computer Science ({LiCS `05})}, pages 311--320. IEEE
  Press, 2005.

\bibitem{Scott1972}
D.~Scott.
\newblock Continuous lattices.
\newblock In {\em Toposes, algebraic geometry and logic}, volume 274 of {\em
  LNM}, pages 97--136. Springer, 1972.

\bibitem{Walters1991}
R.~F.~C. Walters.
\newblock {\em Categories and Computer Science}.
\newblock Carslaw Publications, 1991.

\bibitem{Winskel1995}
G.~Winskel and M.~Nielsen.
\newblock Models for concurrency.
\newblock In S.~Abramsky, D.~Gabbay, and T.~Maibaum, editors, {\em Handbook of
  Logic in Computer Science}, pages 1--148. Oxford University Press, 1995.

\bibitem{Wyler1991}
O.~Wyler.
\newblock {\em Lecture notes on topoi and quasitopoi}.
\newblock World Scientific, 1991.

\end{thebibliography}
\end{document}